\documentclass[10pt,leqno,oneside]{extarticle}

\usepackage{cite}
\usepackage[small]{caption}
\usepackage{amsmath}
\usepackage{amssymb}
\usepackage{amsfonts}
\usepackage{setspace}
\usepackage{latexsym}
\usepackage{mathrsfs}
\usepackage{epsfig}
\usepackage{yfonts}
\usepackage{graphicx}

\newcommand{\n}{\noindent}
\newcommand{\be}{\begin{equation}}
\newcommand{\ee}{\end{equation}}
\newcommand{\ben}{\begin{displaymath}}
\newcommand{\een}{\end{displaymath}}
\newcommand{\bea}{\begin{align}}
\newcommand{\eea}{\end{align}}
\newcommand{\bean}{\begin{align*}}
\newcommand{\eean}{\end{align*}}
\newcommand{\beg}{\begin{gather}}
\newcommand{\eeg}{\end{gather}}
\newcommand{\begn}{\begin{gather*}}
\newcommand{\eegn}{\end{gather*}}
\newcommand{\ep}{\hspace{\stretch{1}}$\blacksquare$}
\newcommand{\et}{\hspace{\stretch{1}}$\blacktriangleleft$}
\newcommand{\etp}{\hspace{\stretch{1}}$\blacktriangleleft\blacksquare$}

\newcommand{\vs}{\vspace{0.2cm}}

\newtheorem{Remark}{Remark}
\newtheorem{Definition}{Definition}
\newtheorem{Proposition}{Proposition}
\newtheorem{Theorem}{Theorem}

\newtheorem{Lemma}{Lemma}

\newcommand{\rt}{\mathbb{R}^{3}}

\newcommand{\Rot}{\mathcal{R}}
\newcommand{\st}{T^{2\blacklozenge}}
\newcommand{\nst}{T^{2\lozenge}}
\newcommand{\torus}{\mathbb{T}^{2}}
\newcommand{\curve}{\mathscr C}

\begin{document}

\n {\huge On Ricci curvature and volume growth in 

\vs
\n dimension three.}

\vspace{0.4cm}
\n {\sc\large Martin Reiris.}

\vs
\n \textsc{\small Max Planck Institute f\"ur Gravitationphysik. \\ Albert Einstein Institut - Germany.}\\e-mail: martin@aei.mpg.de.

\vspace{0.6cm}
\n \begin{minipage}[l]{10cm}
\begin{spacing}{1}
{\small We prove that any complete metric on $\rt$ minus a ball with non-negative Ricci curvature and quadratic Ricci-curvature decay, has cubic volume growth.}
\end{spacing}
\n \end{minipage}

{\fontsize{10.8pt}{12.7pt}\selectfont

\vs
\section{Introduction.}

In the study of non-compact manifolds a simple and at the same time rich invariant worth investigating is 
the rate of volume growth of geodesic spheres. For instance, under some local conditions on the curvature, the rate of volume growth, which is an asymptotic invariant, can provide global information. This is the case for example in Anderson's Gap Theorem \cite{MR1074481}. Other instances where the rate of volume growth plays a relevant role is in the global behavior of harmonic and Green functions \cite{MR2483365} or in the existence and structure of cones at infinity \cite{MR1296356}.

Before stating our main theorem let us recall some terminology. Let $g$ be a complete metric in $\rt$. We say that $g$ has {\it cubic volume growth} if \footnote{Sometimes the terminology cubic {\it volume growth} refers to the condition $\omega_{1}\bar{r}^{3}\leq Vol (B(o,\bar{r}))\leq \omega_{2}\bar{r}^{3}$.} 
\be\label{VCL}
\lim_{\bar{r}\uparrow \infty} \frac{Vol(B(o,\bar{r}))}{\bar{r}^{3}}=\rho>0
\ee
namely the limit exists and is non-zero, where here $o$ is the origin in $\rt$ and $B(o,\bar{r})$ is the geodesic ball of center $o$ and radius $\bar{r}$. If $g$ is a complete metric in $\rt$ minus (say) the unit ball $\mathbb{B}^{3}$, then one can always ``extend" the metric $g$ to all of $\rt$ in such a way that the points in $\partial \mathbb{B}^{3}$ are equidistant to $o$ (\footnote{If $r$ is the usual radial coordinate in $\mathbb{B}^{3}$ then the extension can be written in the form $g=dr^{2} + h(r)$ where $h(r)$ ($r>0$) is a a path of two-metrics on $\mathbb{S}^{2}$ with appropriate values for $h(1),h'(1),h''(1)$ to make the extension $C^{2}$.}). In this sense we say that $g$ (in $\rt\setminus \mathbb{B}^{3}$) has cubic volume growth if the extension to $\rt$ has cubic volume growth. This is the same to say that $\lim Vol({\mathcal T}(\partial \mathbb{B}^{3},\bar{r}))/\bar{r}^{3}$ exists and is non-zero, where ${\mathcal T}(\partial \mathbb{B}^{3},\bar{r})$ is the geodesic tubular-neighborhood (in $\rt\setminus \mathbb{B}^{3}$) of $\partial \mathbb{B}^{3}$ and radius $\bar{r}$. The metric has quadratic curvature decay if $|Ric|\leq \Lambda_{0}/r^{2}$ where $r(p)=dist(p,o)$. The following is the main result of this article.
\begin{Theorem}\label{MT} Let $g$ be a complete metric in $\rt$ minus a ball with non-negative Ricci curvature and quadratic curvature decay. Then $g$ has cubic volume growth.
\end{Theorem} 

A few comments on the hypothesis of the theorem are in order. On $\rt\setminus \mathbb{B}^{3}$ and for $\alpha\in (1/2,1)\cup (1,3/2)$ consider the Riemannian metric
$g=dr^{2}+r^{2\alpha}d\Omega^{2}$, where $d\Omega^{2}$ is the round metric on $\mathbb{S}^{2}$ and $r$ is the standard radial coordinate in $\rt$. Extend $g$ to a spherically symmetric metric in $\rt$.
When $\alpha\in (1/2,1)$ then $|Ric|\leq \Lambda_{0}(\alpha)/r^{2\alpha}$ and $Ric(p)\geq 0$ if $r(p)\geq r_{0}(\alpha)$. Moreover, no matter the value of $\alpha$ in $(1/2,1)$ we have
\ben
\lim_{\bar{r}\uparrow \infty} \frac{Vol(B(o,\bar{r}))}{\bar{r}^{3}}=0
\een
The example shows that the hypothesis ``{\it quadratic curvature decay}" in Theorem \ref{MT} can be hardly weakened (and not removable). 
On the other hand suppose that $\alpha\in (1,3/2)$. Then $|Ric|\leq \Lambda_{0}(\alpha)/r^{2}$ but the Ricci curvature is not non-negative. Moreover, no matter the value of $\alpha$ in $(1,3/2)$ we have
\ben
\lim_{\bar{r}\uparrow \infty} \frac{Vol(B(o,\bar{r}))}{\bar{r}^{3}}=\infty
\een
The example shows that the hypothesis ``{\it non-negative Ricci curvature}" in Theorem \ref{MT} cannot be completely removed. 
In this respect an interesting question is if such hypothesis could be replaced by the much simpler one of ``{\it non-negative scalar curvature}". We point out that examples can be given of complete metrics in $\rt$ with quadratic curvature decay and slow-volume growth, namely with $\rho=0$ (\!\!\cite{MR1755117}).     
Finally, replacing $\rt$ by $\mathbb{R}^{n}$ with $n= 2,4,5,6,\ldots$ and ``cubic" by "Euclidean"(\footnote{In the literature the (unhappy) terminology ``Euclidean volume growth" refers to the condition $\omega_{1} r^{n}\leq Vol(B^{n}_{g}(o,r))\leq \omega_{2} r^{n}$.}) may also make the statement false. For instance the flat product metric on $\mathbb{S}^{1}\times \mathbb{R}^{+}$, which has linear volume growth, shows that it would be false when $n=2$ and the well known Tau-NUT Ricci flat instanton in $\mathbb{R}^{4}$ which has cubic volume growth, shows that it would be false when $n=4$. We do not know at the moment if $n=3$ is the only dimension where the statement holds. 

It is worth mentioning that the relation between volume growth and lower curvature decay has been discussed at least in \cite{MR658471}. Their work however does not overlap with ours, but instead, it complements. This is because \cite{MR658471} argues under the hypothesis $Ric\geq \Lambda_{1}/r^{2}$, $\Lambda_{1}>0$, which turns out to be, if one is working on $\rt\setminus \mathbb{B}^{3}$, incompatible\footnote{For several reasons, for instance $Ric\geq \Lambda_{1}/r^{2}$ implies non-cubic volume growth (see \cite{MR658471} $\S 4$).} with $|Ric|\leq \Lambda_{0}/r^{2}$.       

The idea of the proof, which proceeds by contradiction, is somehow simple. In gross terms one proves that if the volume growth is non-cubic then one can partition $\rt$ into a set of manifolds with a sufficient understanding of their topology to be able to prove that their union is topologically incompatible with $\rt$. All the hypothesis in Theorem \ref{MT}, including the dimension, are strongly used. 
Let us elaborate on the argument a bit more technically and, at the same time, explain the organization of the article. Further explanations have to be found inside the proof and in the main text. After assuming that the volume growth is non-cubic, the proof of Theorem \ref{MT} which starts in pg. \pageref{PMTP} (and ends in pg. \pageref{ETMLF}) goes by writing first $\rt$ as the union of an open set containing the origin and with compact closure, and a set
\be\label{UFIN}
\bigcup_{i=i_{0}}^{i=\infty} M(T^{2o}_{i+1},T^{2o}_{i})
\ee
where every $M(T^{2o}_{i+1},T^{2o}_{i})$ is a compact three-manifold with the tori $T^{2o}_{i+1}$, $T^{2o}_{i}$ as its boundary components (see Figure \ref{UU}). In this union the interiors $M(T^{2o}_{i+1},T^{2o}_{i})^{\circ}$ are pairwise disjoint. The manifolds $M(T^{2o}_{i+1},T^{2o}_{i})$ and the tori $T^{2o}_{i}$, $i\geq i_{0}$ are carefully defined in Section \ref{AD} and Proposition \ref{SADP}. The proof continues by showing that every $M(T^{2o}_{i+1},T^{2o}_{i})$ is an {\it irreducible manifold with incompressible boundary} (shortly IIB-manifold), to conclude finally, due to the topological properties of IIB-manifolds, that (\ref{UFIN}) cannot cover $\rt$ up to an open set of compact closure containing the origin. 
The properties of IIB-manifolds that are required for such conclusion are described in Section \ref{SRTIIB}. In particular it is recalled that any union of IIB-manifolds along boundary components is a IIB-manifold. 
To be concrete, the conclusion, or, more precisely, the contradiction, arises as follows. Pick a coordinate sphere $\mathbb{S}^{2}_{\bar{r}}=\partial \mathbb{B}^{3}(o,\bar{r})$ of coordinate radius $\bar{r}$ in $\rt$, with $\bar{r}$ big enough that $\mathbb{S}^{2}_{\bar{r}}$ is inside the union (\ref{UFIN}) (indeed inside a finite union of $ M(T^{2o}_{i+1},T^{2o}_{i})$'s). As such union is irreducible, the sphere $\mathbb{S}^{2}_{\bar{r}}$ must bound a three-ball in it which forcefully must be $\mathbb{B}^{3}(o,\bar{r})$. Thus the origin must belong to the union (\ref{UFIN}) which is a contradiction.   
That the $M(T^{2o}_{i+1},T^{2o}_{i})$ are IIB-manifolds is deduced during the proof from various informations.
Firstly, the $M(T^{2o}_{i+1},T^{2o}_{i})$ are constructed as finite unions of manifolds $U_{k,l}$ of a special {\it annuli decomposition} ${\mathcal U}$, which, as defined and described in Section \ref{AD}, are particular partitions of $\rt$. The properties of the special annuli decomposition that we will use, are described in detail in Proposition \ref{SADP} in Section \ref{SAD}. 
Roughly speaking, the decomposition is constructed by carefully studying the annuli
\ben
A_{k}(10^{n_{1}},10^{n_{2}}):=B(o,10^{k+n_{2}})\setminus \overline{B(o,10^{k+n_{1}})},\ k=k_{0},k_{0}+2,\ldots
\een	
($n_{1}<n_{2}$ integers but fixed) provided with the scaled metrics $g_{k}:=\frac{1}{10^{2k}}g$ and by means of the Cheeger-Gromov-Fukaya theory of volume collapse with bounded diameter and curvature\footnote{This is possible because over every annulus $A_{k}(10^{n_{1}},10^{n_{2}})$ the $g_{k}$-Ricci curvature is bounded uniformly by $10^{-2n_{1}}\Lambda_{0}$ and because $\lim_{k\uparrow \infty} Vol_{g_{k}}(A_{k}(10^{n_{1}},10^{n_{2}}))=0$, which follows from the fact that, if the volume growth is non-cubic then it has to be sub-cubic as a consequence of the Bishop-Gromov monotonicity 
$Vol({\mathcal T}(\partial \mathbb{B}^{3},\bar{r}))/\bar{r}^{3}\downarrow$, that is $\rho=0$ in (\ref{VCL}). Moreover, as explained in Example I in Section \ref{VCTM}, a certain (but practical) uniform diameter boundednes is obtained for $(A_{k}(10^{n_{1}},10^{n_{2}}),g_{k})$ from Liu's Ball covering property.}. Still such theory for manifolds with boundary has not been appropriately discussed in the literature. To fill in this small gap and to provide a reasonable background for those not familiar with it we dedicate the whole Section \ref{VCBDC} to analyze this matter (see in particular Footnote \ref{FONP} in pg. \pageref{FONP}). 
Secondly, from Proposition \ref{SADP} and crucially Proposition \ref{CRUC} and by using further topological properties of three-manifolds enclosed by embedded tori in $\rt$ which are discussed in Section \ref{SRTIIB}, it is deduced that the pieces $U_{k,l}$ which make up $M(T^{2o}_{i+1},T^{2o}_{i})$ can be grouped appropriately to form IIB-manifolds. Thus it is obtained that every $M(T^{2o}_{i+1},T^{2o}_{i})$ is a union of IIB-manifolds and therefore a IIB-manifold itself.   

We explain in the Appendix a couple of technical propositions whose inclusion inside the text would cause much disruption.  

The article has a good amount of background material, examples and illustrations. 

\section{ Preliminaries.}
\subsection{ Basic notation.}\label{BNS}

\noindent $\mathbb{S}^{n}$, $n\geq 1$ will be the unit sphere in $\mathbb{R}^{n+1}$ and
$\mathbb{T}^{2}=\mathbb{S}^{1}\times \mathbb{S}^{1}$ the two-dimensional torus. $\mathbb{S}^{1}$ and $\mathbb{T}^{2}$ will be thought both as manifolds and as Lie groups.  Furthermore $\mathbb{B}^{n}(o,r)=\{\bar{x}\in \mathbb{R}^{n}, |\bar{x}|<r\}$ will be the open ball of center the origin $o=(0,0,0)$ and radius $r$ ($|\bar{x}|$ is the Euclidean norm of a point $\bar{x}$ of $\mathbb{R}^{n}$). 
$\mathbb{B}^{n}=\mathbb{B}^{n}(o,1)$. $\mathbb{I}=\mathbb{B}^{1}$.

Let $(M,g)$ be a compact connected Riemannian manifold with boundary. The Riemannian metric $g$ induces a metric $d_{g}$ in $M$ (as usual) by
\be\label{MD}
d_{g}(p,q)=dist_{g}(p,q)=\inf\{length_{g}(\gamma(p,q)),\gamma(p,q)\text{ $C^{1}$-curve from $p$ to $q$}\}
\ee
However if  $(\Omega,g)\subset (M,g)$, then on $\Omega$ one can consider two different distances, the distance (\ref{MD}) of $(\Omega,g)$ or the distance (\ref{MD}) of $(M,g)$ restricted to $\Omega$. This situation will appear often and for this reason and to avoid confusion we will denote them by $d^{\Omega}_{g}$ and $d^{M}_{g}$ respectively.

In this article the Riemannian space $(\Omega,g)$ will also denote the metric space $(\Omega,d^{\Omega}_{g})$. 

We will always use the following definitions of {\it diameter} $diam_{g}(\Omega)$ and {\it radius (to the boundary)} $rad_{g}(\Omega)$, even when $(\Omega,g)\subset (M,g)$:
\ben
diam_{g}(\Omega)=\sup\{d^{\Omega}_{g}(p,q);p,q\in \Omega\},\ rad_{g}(\Omega)=\sup\{d^{\Omega}_{g}(p,\partial \Omega); p\in \Omega\}.
\een

Manifold interiors $\Omega\setminus \partial \Omega$ are denoted by $\Omega^{\circ}$. To us a metric ball of center $p\in \Omega^{\circ}$ and radius $r$ is a {\it geodesic ball} if $r<d^{\Omega}_{g}(p,\partial \Omega)$. 

The ends of Theorems, Lemmas or Propositions are marked with $\blacksquare$, while the end of a claim or the end of an Example, is marked with a $\blacktriangleleft$.

\subsection{ Surfaces in $\mathbb{R}^{3}$ and irreducible three-manifolds with incompressible boundary.}\label{SRTIIB}
From now on we let $S$ be a smoothly embedded compact, orientable and boundaryless surface in $\mathbb{R}^{3}$. Any $S$ divides $\rt$ into two open connected components.  
We will denote by $M(S)$ the closure of the bounded component. For instance if $S\sim \mathbb{S}^{2}$ then $S$ bounds a three-ball \cite{Alexander1}.   
If $S\cap S'=\emptyset$ then either
\be\label{TO}
M(S)\cap M(S')=\emptyset,\ M(S)\subset M(S')^{\circ}, \text{ \   or }\ M(S')\subset M(S)^{\circ}.
\ee
Moreover if $S'\subset M(S)^{\circ}$ then $S$ belongs to $\rt\setminus M(S')$ and therefore $M(S')\subset M(S)^{\circ}$. In particular if $S'\sim \mathbb{S}^{2}$ and $S'\subset M(S)^{\circ}$ then $S'$ bounds a three-ball inside $M(S)^{\circ}$. Recall that a three-manifold is {\it irreducible} if every embedded two-sphere bounds a three-ball. Thus for any $S$, $M(S)$ is an irreducible manifold. 

If $S\sim \mathbb{T}^{2}$ then either $M(S)$ is a solid torus, i.e. $\sim \bar{\mathbb{B}}^{2}\times \mathbb{S}^{1}$, or, $S=\partial M(S)$ is incompressible in $M(S)$, where recall, $N$ is an {\it incompressible} boundary component of a manifold $M$ if $i_{*}:\pi_{1}(N)\rightarrow \pi_{1}(M)$ is injective (here $i:N\rightarrow M$ is the inclusion). To see this think $S$ as a surface in $\mathbb{S}^{3}$ via $S\subset \rt \subset \big(\rt \cup \{\infty\}\big)\sim \mathbb{S}^{3}$.  If $M(S)$ is a solid torus we are done. If not then $\mathbb{S}^{3}\setminus M(S)^{\circ}$ is a solid torus (this is due to Alexander \cite{Alexander1}). If $\mathbb{S}^{3}\setminus M(S)^{\circ}$ represents the unknot then $M(S)$ is a solid torus but we are assuming that it is not. Then $\mathbb{S}^{3}\setminus M(S)^{\circ}$ is not the unknot. Theorem 11.2 in \cite{MR1472978} shows that in this case $S$ is incompressible in $M(S)$ as claimed. Summarizing, for any $S\sim \mathbb{T}^{2}$, $M(S)$ is either a solid torus or an irreducible manifold with incompressible boundary.  

Other examples of irreducible manifolds with incompressible boundary components (in short, ``IIB" manifolds) are compact {\it Seifert} manifolds with at least two boundary components (\!\!\cite{MR705527} pgs. 431-432 and Corollary 3.3). Recall that 
a Seifert manifold is one admitting a foliation ${\mathcal C}$ by circles $\mathscr C$ around any of which there is a fibered neighborhood isomorphic to a {\it fibered solid torus or Klein bottle} (see \cite{MR705527}, pg. 428). The class of IIB manifolds is closed under sums along boundary components. Precisely we have (Lemma 1.1.4 in \cite{MR0224099})
\begin{Proposition}\label{IIBM}\
\begin{enumerate}
\item[\bf I.]  Let $M_{1}$ and $M_{2}$ be two {\rm IIB} manifolds and let $f:N_{1}\rightarrow N_{2}$ be a diffeomorphism between a boundary component $N_{1}$ of $M_{1}$ and a boundary component $N_{2}$ of $M_{2}$. Then the manifold which results from identifying through $f$ the boundary $N_{1}$ of $M_{1}$ to the boundary $N_{2}$ of $M_{2}$ is a {\rm IIB} manifold. 
\item[\bf II.] Let $M_{1}$ be a {\rm IIB} manifold and let $f:N_{1}\rightarrow N_{2}$ be a diffeomorphism between the boundary components $N_{1}\neq N_{2}$ of  $M_{1}$. Then, the manifold which results from identifying through $f$ the boundary $N_{1}$ to the boundary $N_{2}$ of the manifold $M_{1}$ is a {\rm IIB} manifold.
\end{enumerate}
Therefore, any sum of {\rm IIB} manifolds along any number of boundary components is a {\rm IIB} manifold.
\end{Proposition}
Yet, there is a simple but important situation when the sum of a IIB manifold and a non-IIB manifold results in a IIB manifold. The case is when $M_{1}$ is a Seifert manifold with Seifert structure ${\mathcal C}$ and at least three-boundary components, $M_{2}$ is a solid torus and the gluing function $f:N_{1}(\subset \partial M_{1})\rightarrow N_{2}(=\partial M_{2})$ send circles ${\mathscr C}$ in ${\mathcal C}$ into non-contractible circles $f({\mathscr C})$ as a circle in $M_{2}$. The reason is that in this situation the $\mathbb{S}^{1}$-foliation $f({\mathcal C})$ of $N_{2}=\partial M_{2}$ can always be extended to a Seifert structure in $M_{2}$ and thus making the sum a Seifert manifold with at least two boundary components and therefore a IIB manifold. To construct the extension of $f({\mathcal C})$ proceed as follows. On $M_{2}\sim \mathbb{B}^{2}\times \mathbb{S}^{1}$ denote points by $(\bar{x},s)$, $\bar{x}\in \mathbb{B}^{2}$ and $s\in \mathbb{S}^{1}$. Then, for any $r\in [0,1]$ define $F_{r}:\mathbb{B}^{2}\times \mathbb{S}^{1}\rightarrow \mathbb{B}^{2}\times \mathbb{S}^{1}$ by $F_{r}(\bar{x},s)=(r\bar{x},s)$. The desired extension of $f({\mathcal C})$ is $\{F_{r}(\curve);\curve\in {\mathcal C},r\in [0,1]\}$.
\subsection{ Annuli decompositions.}\label{AD}
Let $g$ be a complete metric in $\rt$. For every $b>a>0$ we let $A_{g}(a,b)=B_{g}(o,b)\setminus \overline{B_{g}(o,a)}$ be the {\it (open) annulus} with radii $a$ and $b$ and center the origin $o$. 
\begin{Definition}\label{D1} A set ${\mathcal{U}}=\{U_{k,l}; k=k_{0}+2j, j=0,1,2,3,..., l=1,2,\ldots,l(k) \}$
of compact three-submanifolds of $\rt$ with smooth boundary is an annuli decomposition iff the following conditions are fulfilled:  
\begin{enumerate}
\item $U_{k,l}\subset A_{g}(10^{k-1},10^{k+3})$,
\item $U_{k,l}\cap A_{g}(10^{k-1},10^{k})\neq \emptyset$, and $U_{k,l}\cap A_{g}(10^{k+1},10^{k+2})\neq \emptyset$, 
\item $\partial U_{k,l}\subset \big(A_{g}(10^{k-1},10^{k})\cup A_{g}(10^{k+1},10^{k+2})\big)$,
\item $\text{If } (k,l)\neq (k',l')$, then $U_{k,l}^{\circ}$ and $U_{k',l'}^{\circ}$ are disjoint and if $U_{k,l}$ and $U_{k',l'}$ intersect then they do in a set of boundary components (of both, $U_{k,l}$ and $U_{k',l'}$),
\item $U_{k_{0}-2}:=\rt\setminus \big( \bigcup \limits_{U_{k,l}\in {\mathcal U}} U_{k,l}\big)^{\circ}$ is compact.
\end{enumerate}
\end{Definition}
\begin{figure}[h]
\centering
\includegraphics[width=13cm,height=10cm]{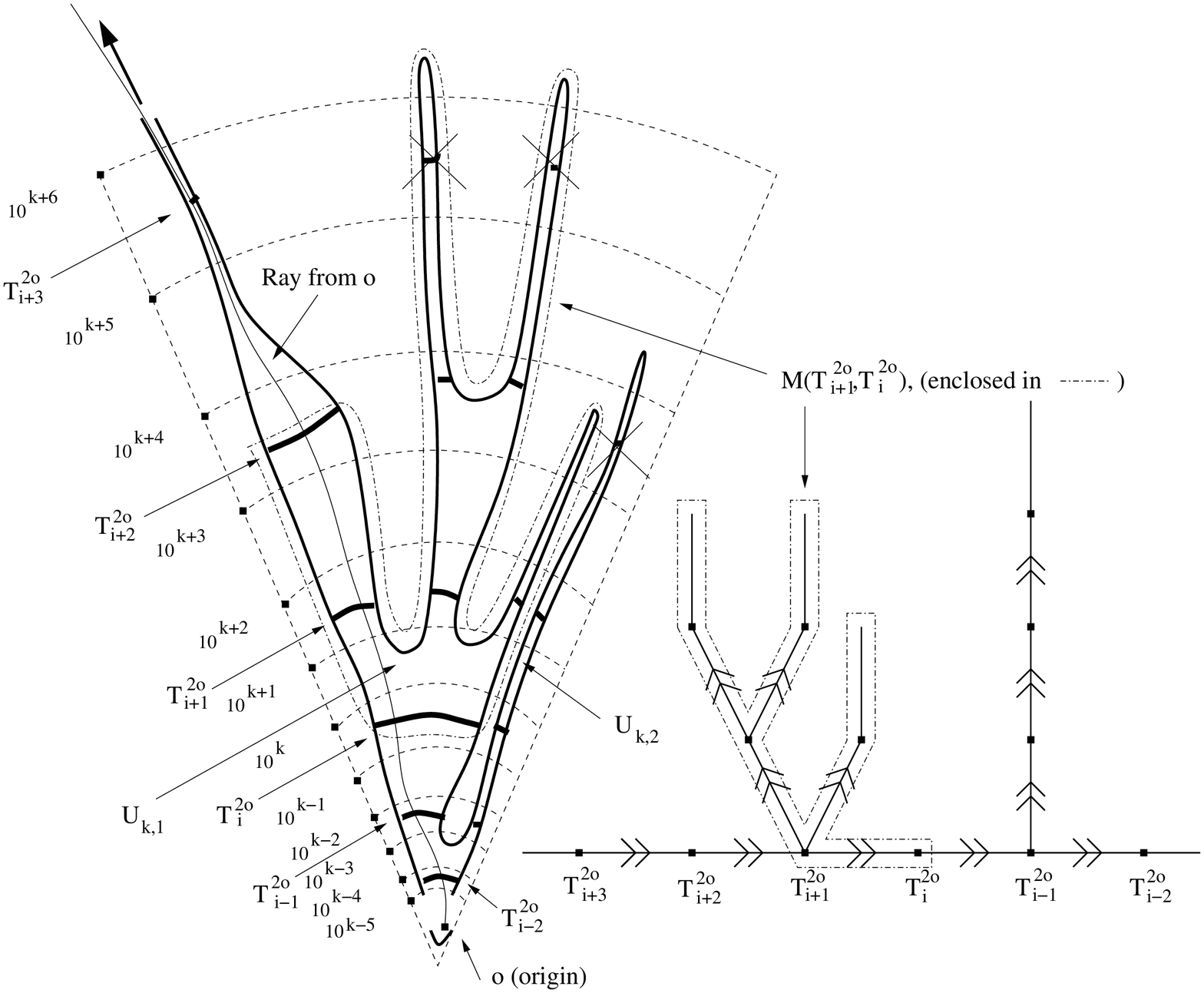}
\caption{On the left side the figure schematizes a part of an annuli decomposition. We have indicated only the pieces $U_{k,1}$ and $U_{k,2}$ but but every region enclosed by thick lines is a piece $U_{k,l}$. We have also explicitly indicated the surfaces $S^{o}_{i}$, however with a $T^{2o}_{i}$, as this is the notation to be used in Section \ref{PMT} where the proof of Theorem \ref{MT} is carried out.  On the right side is represented the corresponding part of the tree induced by the order $\ll$. On both sides we have enclosed in a dash/point line ($-\cdot-\cdot-$) the manifold $M(T^{2o}_{i+1},M(T^{2o}_{i}))$. On the left it can also be seen crossed thick lines. This is for a later exemplification when in Proposition \ref{SADP} in Section \ref{SAD} we explain how the special annuli decomposition, to be used in the proof of Theorem \ref{MT}, is constructed. The cross indicates that such ``cuts", as we refer them there, are to be discarded.}
\label{UU}
\end{figure} 
Let ${\mathcal N}$ be the set of boundary components of the manifolds $U_{k,l}$ in an annuli decomposition ${\mathcal U}$. Elements of ${\mathcal N}$ are pairwise disjoint compact, orientable and embedded surfaces. We can order them as follows: $S\ll S'$ iff $M(S)\subset M(S')$. The order is not necessarily a linear order, as there can be two elements not related. However there is an important subset which is linearly ordered, this is the set
${\mathcal N}^{o}=\{S\in {\mathcal N}; o\in M(S)\}$ (use (\ref{TO})).  
Thus ${\mathcal N}^{o}=\{S^{o}_{1},S^{o}_{2},S^{o}_{3},\ldots\}$ with $S^{o}_{1}\ll S^{o}_{2}\ll S^{o}_{3}\ll \ldots$. We will be using this notation (the upper-index $o$ is from ``origin"). We will also use later the notation $M(S^{o}_{i},S^{o}_{i'}):=M(S^{o}_{i})\setminus M(S^{o}_{i'})^{\circ}$ for the region enclosed by $S^{o}_{i}$ and $S^{o}_{i'}$, $i>i'$.  
Moreover because of (\ref{TO}) we have the following property: given two elements $S\ll S'$ in ${\mathcal N}$ then there is a unique (and finite) maximal ``chain" $\{S_{0},\ldots, S_{n}\}\subset {\mathcal N}$ such that $S=S_{0}\ll S_{1}\ll\ldots \ll S_{n-1}\ll S_{n}=S'$. Later, in the proof of Theorem \ref{MT}, we will use the notation $\{S,S'\}\rightarrow \{S,S_{1},\ldots,S_{n-1},S'\}$. 

A representation of an annuli decomposition can be seen in Figure \ref{UU}. The figure shows also the tree induced by the order $\ll$. 

We note in passing that the notion of annuli decomposition (see also the notion of $(\underline{\epsilon},\overline{\epsilon})$-connected components in Definition \ref{ECCD}) and that of ``{\it chopping}" defined in \cite{MR1173034} share some similarities. 
\subsection{ Volume collapse with bounded diameter and curvature.}\label{VCBDC}
\subsubsection{ The Gromov-Hausdorff distance and a relevant example.}\label{GHDRE}
The {\it Gromov-Hausdorff distance} (shortly GH-distance) \cite{MR2307192} between two compact metric spaces $(X,d_{X})$ and $(Y,d_{Y})$ is defined as the infimum of the $\delta>0$ such that there exists, on the disjoint union $X\sqcup Y$, a metric $d_{X\sqcup Y}$ extending $d_{X}$ and $d_{Y}$ such that 
\be\label{GHI}
Y\subset {\mathcal T}_{d_{X\sqcup Y}}(X,\delta)\ \text{and } X\subset {\mathcal T}_{d_{X\sqcup Y}}(Y,\delta)
\ee
where ${\mathcal T}_{d_{X\sqcup Y}}(X,\delta)$ and ${\mathcal T}_{d_{X\sqcup Y}}(Y,\delta)$ are the $d_{X\sqcup Y}$-metric neighborhoods of $X$ and $Y$ and radius $\delta$, respectively. 

We introduce now some terminology to be used during the rest of the article. We will say that a sequence of compact manifolds $(M_{i},g_{i})$ {\it metrically collapses} to a space $(X,d)$ if it converges in the GH-topology to $(X,d)$ and the Hausdorff dimension of $X$ is less than that of $M_{i}$ (which we assume is constant). 
If the GH-distance between $(M,g)$ and $(X,d)$ is less or equal than $\epsilon$ then we say that $(M,g)$ is $\epsilon${\it -close to} $(X,d)$. If the GH-distance between $(M,g)$  and a point is less or equal than $\epsilon$ we say that $(M,g)$ is $\epsilon${\it-collapsed} (for the distance of $(M,g)$ to a point see \cite{MR2243772}). 
   
We present below an example where we estimate the distance between two metric spaces that will be relevant to us in the proof of the Step {\bf C} inside the proof of the Proposition \ref{CRUC}. 

\vs
\n {\bf Example of a Gromov-Hausdorff distance estimation.} Let $I$ be a compact interval in $\mathbb{R}$ of length $|I|\geq 1$. Let $h$ be a flat metric in $\mathbb{T}^{2}$ of diameter $\Gamma$. Provide $X=\mathbb{T}^{2}\times I$ with the metric $d_{X}$ induced from the Riemanian flat product-metric $g=dx^{2}+h$. Intuitively, if $\Gamma$ is small then $(X,d_{X})$ should be close metrically to the interval $I$. More precisely, it should be close to the metric space $(Y,d_{Y})=(I,|\ |)$ where $d_{Y}(x_{1},x_{2})=|x_{1}-x_{2}|$. We show now the following upper and lower bounds for the GH-distance between $(X,d_{X})$ and $(Y,d_{Y})$,  (when $\Gamma\leq 1$)
\be\label{DGHD}
\frac{\Gamma}{5}\leq dist_{GH}(X,Y)\leq \frac{\Gamma}{2}
\ee

$\bullet$ {\it The upper bound.} Points in $\mathbb{T}^{2}$ are denoted by $t$, points in $I$ by $x$, and thus points in $X=\mathbb{T}^{2}\times I$ by $(t,x)$. Let $t_{0}$ be a point in $\mathbb{T}^{2}$ such that $\overline{B_{h}(t_{0},\Gamma/2)}=\mathbb{T}^{2}$ (such point always exists). 
For every $\epsilon>0$ define the distance $d^{\epsilon}_{X\sqcup Y}$ as equal to $d_{X}$ and $d_{Y}$ when restricted to $X$ and $Y$ respectively and as $d^{\epsilon}_{X\sqcup Y}((t,x),x')=d_{X}((t,x),(t_{0},x'))+\epsilon$ for the distance between $(t,x)\in \mathbb{T}^{2}\times I$ and $x'\in I$. Now, (\ref{GHI}) holds for $\delta(\epsilon)=\Gamma/2+2\epsilon$ and for any $\epsilon>0$. Therefore $dist_{GH}(X,Y)\leq \frac{\Gamma}{2}$. 

$\bullet$ {\it The lower bound.} Make $dist_{GH}\big(X,Y)\big)=\Gamma/\mu$ for a $\mu$ that we will estimate as $\mu<5$. Let $t_{1}$ and $t_{2}$ be two points in $\mathbb{T}^{2}$ such that $dist_{h}(t_{1},t_{2})=\Gamma$. Let also $p_{1}=(t_{1},0),\ p_{2}=(t_{2},0),\ p_{3}=(t_{1},\Gamma),\ p_{4}=(t_{2},\Gamma)$, forming an ``square" in $X$: i.e. 
$d_{X}(p_{1},p_{2})=d_{X}(p_{2},p_{4})=d_{X}(p_{4},p_{3})=d_{X}(p_{3},p_{1})=\Gamma$ and $
d_{X}(p_{1},p_{4})=d_{X}(p_{2},p_{3})=\sqrt{2} \Gamma$.
By the definition of the GH-distance, for every $\epsilon>0$ there is $d^{\epsilon}_{X\sqcup Y}$ extending $d_{X}$ and $d_{Y}$, and satisfying (\ref{GHI}) with $\delta(\epsilon)=\Gamma/\mu+\epsilon$. Therefore there are points $x_{1},x_{2},x_{3}$ and $x_{4}$ in $I$ such that for every $j=1,2,3,4$ we have
$d^{\epsilon}_{X\sqcup Y}(p_{j},x_{j})\leq \frac{\Gamma}{\mu}+\epsilon$. From this and the triangle inequalities
\begin{gather*}
d_{Y}(x_{j},x_{k})\leq d^{\epsilon}_{X\sqcup Y}(x_{j},p_{j})+d_{X}(p_{j},p_{k})+d^{\epsilon}_{X\sqcup Y}(p_{k},x_{k}),\\
d_{X}(p_{j},p_{k})\leq d^{\epsilon}_{X\sqcup Y}(x_{j},p_{j})+d_{Y}(x_{j},x_{k})+d^{\epsilon}_{X\sqcup Y}(p_{k},x_{k})
\end{gather*}
we get, when  $(j,k)$ is not $(1,4)$ or $(2,3)$ 
\be
\label{GH1}\ |x_{j}-x_{k}|\leq 2\frac{\Gamma}{\mu}+\Gamma+2\epsilon, \text{ and } \Gamma\leq 2\frac{\Gamma}{\mu} + |x_{j}-x_{k}|+2\epsilon,
\ee
while when $(j,k)$ is $(1,4)$ or $(2,3)$
\be
\label{GH2} 
\sqrt{2} |x_{j}-x_{k}|\leq 2\frac{\Gamma}{\mu}+\sqrt{2} \Gamma+2\epsilon, \text{ and } \sqrt{2}\Gamma\leq 2\frac{\Gamma}{\mu} + |x_{j}-x_{k}|+2\epsilon,
\ee
We will use inequalities (\ref{GH1}) and (\ref{GH2}) in what follows. 
Suppose that $x_{1}\leq x_{3}$ (the case $x_{1}\geq x_{3}$ is symmetric). Then:

\ \ \ - If $x_{4}\leq x_{3}$ we have $|x_{1}-x_{4}|=||x_{1}-x_{3}|-|x_{3}-x_{4}||$ which using (\ref{GH1}) is less or equal than $4\Gamma/\mu+4\epsilon$, i.e. $|x_{1}-x_{4}|\leq 4\Gamma/\mu+4\epsilon$. On the other hand from this and (\ref{GH2}) we obtain $\sqrt{2} \Gamma\leq 6\Gamma/\mu+6\epsilon$, for every $\epsilon>0$ and therefore $\mu<5$. 

\ \ \ - If $x_{4}\geq x_{3}$ then we have two possibilities (i) $x_{2}\geq x_{1}$ or (ii) $x_{2}\leq x_{1}$. (i) is symmetric to the case we have considered before under the change $x_{1}\rightarrow x_{3}$, $x_{3}\rightarrow x_{1}$ and $x_{4}\rightarrow x_{2}$. We consider then (ii). In this case we have $|x_{2}-x_{4}|= |x_{2}-x_{1}|+|x_{1}-x_{3}|+|x_{3}-x_{4}|$ which by (\ref{GH1}) is greater or equal than $3\Gamma-6\Gamma/\mu - 6\epsilon$, i.e. $|x_{2}-x_{4}|\geq 3\Gamma-6\Gamma/\mu-6\epsilon$. From this and (\ref{GH1}) again we obtain $4\Gamma/\mu\geq \Gamma-4\epsilon$, for every $\epsilon>0$ and therefore $\mu<5$. \et
\subsubsection{ The local models of collapse and examples.}
Locally there are only five types of models describing the metric limit of boundaryless compact three-manifolds collapsing in volume with curvature and diameter bounds. If $(X,d)$ is a limit metric space and $x\in X$ then either $x$ is the only point of $X$ or there is a neighborhood of $x$ locally isometric to one of the following four possibilities:
\begin{enumerate}
\item[{\bf I.a\ }] an interval $I=(-a,a)$, with $-a<x=0<a$, provided with the standard metric $d(x_{1},x_{2})=|x_{1}-x_{2}|$,
\item[{\bf I.b\ }] an interval $I=[0,a)$, with $x=0<a$, provided with the standard metric $d(x_{1},x_{2})=|x_{1}-x_{2}|$,
\item[{\bf II.a}] a disc $D=\mathbb{B}^{2}(o,a)$, $x=o$, provided with a metric $d$ induced from a $C^{1,\beta}$-Riemannian metric,
\item[{\bf II.b}] a disc $D=\mathbb{B}^{2}(o,a)$, $x=o$, provided with a metric $d$ induced from the quotient of a $C^{1,\beta}$-Riemannian metric by the action of $\mathbb{Z}_{q}$, $q\geq 1$ by isometries leaving the origin $o$ fixed. 
\end{enumerate}
The point $x=0$ in case ${\bf I.a}$ and the point $x=o$ in case ${\bf II.b.}$ will be here called {\it singular points} and denoted by $Sing(X)$.  A manifold locally of the form {\bf II.a} or {\bf II.b} will be called a $C^{1,\beta}$ orbifold.
\begin{figure}[h]
\centering\includegraphics[width=9cm,height=1.7cm]{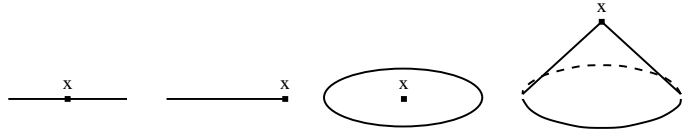}\caption{The local models of collapse. From left to right, models: {\bf I.a}, {\bf I.b}, {\bf II.a}, {\bf II.b}.}\label{LM}
\end{figure} 

That {\bf I.a, I. b, II.a,} and {\bf II.b} are the only possible models is an important consequence of the Cheeger-Gromov-Fukaya theory of collapse under curvature bounds \cite{MR950552}. We comment on this in what follows. First, the limit space is always of  integer dimension and therefore if it not a point it must be of dimension one or two as stated in Theorem 0.6 in pg. 2 (and the paragraph below it) of \cite{MR950552}. That when the dimension is two the models are of the forms {\bf II.a} and {\bf II.b} is the content of Proposition 11.5 in pg. 186 in \cite{MR1145256} (Proposition 11.5 is a Corollary to Theorem 11.1 in pg. 184, which is a restatement of Theorem 0.6 in \cite{MR950552}). That when the dimension is one the models are of the forms {\bf I.a} and {\bf I.b} follows from Theorem 0.5 in pg. 2 of \cite{MR950552} after Definition 0.4. Indeed by Theorem 0.5 and Definition 0.4 there is a neighborhood of $x$ homeomorphic to the quotient of $\mathbb{B}^{m}$ (with $x=o$ and for some $m$) by the action of a Lie subgroup of $O(3)$. Thus, if the Hausdorff dimension is one then the space must be of the types {\bf I.a} or {\bf I.b}, as these are the only possible quotients of dimension one. Note that it is excluded for instance the union of three or more segments in a point (if we remove $x=o$ from the quotient the space must be still connected).     

Below we are going to give four examples showing how everyone of the four cases above can be realized. They are illustrative and do not play any other role in the article. For this reason the presentation is rather synthetic. The examples give sequences of Riemannian manifolds $(M_{n},g_{n})$ converging to a $(X,d)$ as in {\bf I.a}, {\bf I.b}, {\bf II.a} and {\bf II.b} (in this order). We define first the sequence $(M_{n},g_{n})$ and give what is going to be the limit space $(X,d)$. After the definitions for every one of the cases {\bf I.a, I.b, II.a} and {\bf II.b} are made, we list the  geometric properties of the convergence process applying to each. The justifications are just computational and because they play no role in the article are left to the readers. Finally let us mention that the examples show essentially all what can occur locally in volume collapse with curvature and diameter bounds besides collapse to a point (see Lemma \ref{LL}).

We will use the following notation. 
The rotational group of $\mathbb{R}^{2}\sim \mathbb{C}$ will be denoted by $\Rot$. Obviously $\mathbb{U}(1)\sim\Rot$ under the homomorphism $u\in \mathbb{U}(1)\rightarrow R(u)\in \Rot$, with $R(u)z=uz$ for any $z\in \mathbb{C}$.
Also for any natural number $q$ let $\Rot^{q}\sim \mathbb{Z}_{q}$ be the subgroup of rotations generated by $R(e^{2\pi i/q})$.
Finally the group of rotations on the first factor $\mathbb{R}^{2}$ in $\mathbb{R}^{2}\times \mathbb{R}^{2}$ with be denoted by $\Rot_{1}$ and the group of rotations on the second factor will be denoted by $\Rot_{2}$.
Note that the set $\mathbb{B}^{2}\times \mathbb{S}^{1}\subset \mathbb{R}^{2}\times \mathbb{R}^{2}$ and the set $\torus \subset \mathbb{R}^{2}\times \mathbb{R}^{2}$ are invariant under $\Rot_{1}\times \Rot_{2}$. In particular $\torus\times \mathbb{I}\subset \mathbb{R}^{2}\times \mathbb{R}^{2}\times \mathbb{R}$ is invariant under $\Rot_{1}\times \Rot_{2}$.

\vs
\n {\bf Example I.a.} 

\vspace{0.1cm}
\n $\bullet$ $(M_{n},g_{n})$ - Let $\tilde{M}=\torus\times \mathbb{I}$ and provided with a smooth and $\Rot_{1}\times \Rot_{2}$-invariant Riemannian metric $\tilde{g}$. Let $G_{n}\sim \mathbb{Z}_{n}\times \mathbb{Z}_{n}$ be the group generated by the rotations $\Rot_{1}(e^{2\pi i/n})$, $\Rot_{2}(e^{2\pi i/n})$. Let $M_{n}=\tilde{M}/G_{n}$ be the quotient of $\tilde{M}$ by $G_{n}$, $\pi_{n}:\tilde{M}\rightarrow M_{n}$ the covering map and $g_{n}$ the projected metric on $M_{n}$, namely $\pi^{*}_{n}(g_{n})=\tilde{g}$. 

\vspace{0.1cm}
\n $\bullet$ $(X,d)$ - Let $X=\torus\times \mathbb{I}/(\Rot_{1}\times \Rot_{2})$ with the induced quotient metric $d$ and let $f_{n}:M_{n}\rightarrow X$ be the projection.  

\vs
\n {\bf Example I.b.}

\vspace{0.1cm}
\n $\bullet$ $(M_{n},g_{n})$ - Let $\tilde{M}=\mathbb{B}^{2}\times \mathbb{S}^{1}$ and provided with a smooth and $\Rot_{1}\times \Rot_{2}$-invariant Riemannian metric $\tilde{g}$. Let $G_{n}\sim \mathbb{Z}_{n^{2}}$ be the group generated by the rotations $\Rot_{1}(e^{2\pi i/n})\times \Rot_{2}(e^{2\pi i/n^{2}})$. Let $M_{n}=\tilde{M}/G_{n}$ be the quotient of $\tilde{M}$ by $G_{n}$, $\pi_{n}:\tilde{M}\rightarrow M_{n}$ the covering map and $g_{n}$ the projected metric on $M_{n}$, namely $\pi^{*}_{n}(g_{n})=\tilde{g}$. 

\vspace{0.1cm}
\n $\bullet$ $(X,d)$ - Let $X=\mathbb{B}^{2}/(\Rot_{1}\times \Rot_{2})$ with the induced quotient metric $d$ and let $f_{n}:M_{n}\rightarrow X$ be the projection.  

\vs
\n {\bf Example II.a.} 

\vspace{0.1cm}
\n $\bullet$ $(M_{n},g_{n})$ - Let $\tilde{M}=\mathbb{B}^{2}\times \mathbb{S}^{1}$ and provided with a smooth and $\Rot_{2}$-invariant Riemannian metric $\tilde{g}$. Let $G_{n}\sim \mathbb{Z}_{n}$ be the subgroup of $\Rot_{2}$ generated by the rotations $\Rot_{2}(e^{2\pi i/n})$. Let $M_{n}=\tilde{M}/G_{n}$ be the quotient of $\tilde{M}$ by $G_{n}$, $\pi_{n}:\tilde{M}\rightarrow M_{n}$ the covering map and $g_{n}$ the projected metric on $M_{n}$, namely $\pi^{*}_{n}(g_{n})=\tilde{g}$. 

\vspace{0.1cm}
\n $\bullet$ $(X,d)$ - Let $X=\mathbb{B}^{2}$, with the induced quotient metric $d$. Let $f_{n}:M_{n}\rightarrow X$ be the projection.  

\vs
\n {\bf Example II.b.} 

\vspace{0.1cm}
\n $\bullet$ $(M_{n},g_{n})$ - Let $\tilde{M}=\mathbb{B}^{2}\times \mathbb{S}^{1}$ provided with a smooth and $\Rot_{1}\times \Rot_{2}$-invariant Riemannian metric $\tilde{g}$. Let $0<p<q$ be two relatively prime natural numbers and let $G_{n}\sim \mathbb{Z}_{qn}$ be the subgroup of $\Rot_{1}\times \Rot_{2}$ generated by the rotations $\Rot_{1}(e^{2\pi p i/q}) \times  \Rot_{2}(e^{2\pi i/qn})$. Let $M_{n}=\tilde{M}/G_{n}$ be the quotient of $\tilde{M}$ by $G_{n}$, $\pi_{n}:\tilde{M}\rightarrow M_{n}$ the covering map and $g_{n}$ the projected metric on $M_{n}$, namely $\pi^{*}_{n}(g_{n})=\tilde{g}$. 

\vspace{0.1cm}
\n $\bullet$ $(X,d)$ - Let $X=\mathbb{B}^{2}/\Rot^{q}$, with the induced quotient metric $d$. Let $f_{n}:M_{n}\rightarrow X$ be the projection.  

\vs
\vs
\n With these definitions for the examples {\bf I.a}, {\bf I.b}, {\bf II.a} and {\bf II.b} it is straightforward to check that,

\begin{enumerate}
\item $Sing(X)=\emptyset$ in cases {\bf I.a}, {\bf II.a} and $Sing(X)=\{o\}$ in cases {\bf I.b} and {\bf II.b}. 
\item In every example the sequence $(M_{n},g_{n})$ converges in the GH-topology to $(X,d)$. The group $G_{n}$ of Deck transformations on $\tilde{M}$ converges to $G=\Rot_{1}\times \Rot_{2}\sim \torus$ in cases {\bf I.a} and {\bf I.b}, to $G={\mathcal R}_{2}\sim \mathbb{S}^{1}$ in case {\bf II.a} and to $G:=\Rot_{1}^{q}\times \Rot_{2}\sim \mathbb{Z}_{q}\times \mathbb{S}^(1)$ in case {\bf II.b}. Moreover $X=\tilde{M}/G$. Let $\pi:\tilde{M}\rightarrow X$ be the projection. Then $Centr_{G}(\pi^{-1}(Sing(X)))=\Rot^{q}_{1}$ where $Centr$ is the centralizer.     
\item In every example $f_{n}:M_{n}\rightarrow X$ is a fibration and $length_{g_{n}}(f_{n}^{-1}(x))\rightarrow 0$. Moreover $f_{n}:M_{n}\setminus f_{n}^{-1}(Sing(X))\rightarrow X\setminus Sing(X)$ is a $\torus$-fiber bundle in cases {\bf I.a, I.b} and a $\mathbb{S}^{1}$-fiber bundle in cases {\bf II.a} and {\bf II.b}. $Centr(Sing(X))$ acts freely on $f_{n}^{-1}(x)$ for any $x\in X\setminus Sing(X)$ and $f_{n}^{-1}(Sing(X))\sim f^{-1}_{n}(x)/Centr(\pi^{-1}(Sing(X)))$. $\blacktriangleleft$
\end{enumerate}
\subsubsection{ Volume collapse of three-manifolds with boundary and with curvature and diameter bounds - an statement.}\label{VCTM}
We discuss now briefly what we will mean by {\it three-manifolds with non-necessarily smooth boundary}. The reader should keep in mind that the notion is just for the purpose of working later with some necessary generality, with no intention whatsoever in developing a new concept, which, as a matter of fact, would be here purposeless. Let $M$ be a compact set on an open manifold $P$. Then we say that $M$ is a {\it compact manifold with non-necessarily smooth boundary} (shortly, {\it manifold with NNSB}) if $M$ is equal to the closure (in $P$) of its interior (in $P$). In this sense the {\it boundary} $\partial M$ of $M$ is defined as $M$ minus the topological interior of $M$ (in $P$) and the {\it manifold's interior} $M^{\circ}:=M\setminus \partial M$ therefore coincides with the topological interior (in $P$). Note that we do not assume that $M^{\circ}$ is connected. 
A subset of $M$ is {\it a submanifold with NNSB} if it is a manifold with NNSB as a subset of $P$. Of course any compact manifold with smooth boundary is a manifold with NNSB. If $P$ carries a Riemannian metric $g$ then we say that $(M,g)$ is a Riemannian manifold with NNSB. In this case the Riemannian metric $g$ induces a metric $d=d^{M}_{g}$ in every connected component of $M^{\circ}$. 
For the discussion below we do not need to extend $d$ to a metric on $M^{\circ}$. 
The distance from a point $p\in M^{\circ}$ to $\partial M$ can be defined in several equivalent and natural ways. For instance $d(p,\partial M)$ as the supremum of the radius of the geodesic balls of center $p$, lying entirely in $M^{\circ}$. Then $d(p,\partial M)$ is realized by the $g$-length of a geodesic starting at $p$, ending at $\partial M$ and whose interior lies in $M^{\circ}$. Define the tubular neighborhoods ${\mathcal T}_{d}(\partial M,\epsilon):=\partial M\cup \{p\in M^{\circ}, d(p,\partial M)<\epsilon\}$.     

\begin{Definition}\label{DEF2} Let ${\mathfrak N}_{0}:\mathbb{R}^{+}\times \mathbb{R}^{+}\rightarrow \mathbb{R}^{+}$ be a non-necessarily continuous function. Then define ${\mathcal M}({\mathfrak{N}}_{0})$ as the set of compact Riemannian manifolds with NNSB $(M,g)$, such that for any $1>\epsilon_{0}>2\epsilon_{1}>0$, the minimum number of geodesic balls of radius $\epsilon_{1}$ covering $M\setminus {\mathcal T}_{d}(\partial M,\epsilon_{0})$ is bounded above by ${\mathfrak N}_{0}(\epsilon_{0},\epsilon_{1})$.
\end{Definition} 

\begin{Remark}\label{R1} The values of ${\mathfrak N}_{0}$ outside the set $\{(\epsilon_{0},\epsilon_{1}), 1>\epsilon_{0}>2\epsilon_{1}>0\}$ are of no relevance. 
\end{Remark}

We would like to comment briefly about the reason of this definition. Recall that given $\Lambda_{0}>0$, $D_{0}>0$ there is ${\mathfrak N}_{0}:\mathbb{R}^{+}\rightarrow \mathbb{R}^{+}$, depending on them, such that for any compact boundaryless Riemannian three-manifold with $|Ric|\leq \Lambda_{0}$, $diam_{g}(M)\leq D_{0}$ the minimum number of balls of radius $\epsilon$ covering $M$ is bounded above by ${\mathfrak N}_{0}(\epsilon)$ (this is due to Gromov; see \cite{MR2243772}, pg. 281). Moreover the existence of such ${\mathfrak N}_{0}$ is equivalent to the precompactness of the family of compact and boundaryless Riemmanian three-manifolds with $|Ric|\leq \Lambda_{0}$ and $diam_{g}(M)\leq D_{0}$, as a set inside the family of compact metric spaces provided with the GH-topology (\!\!\cite{MR2243772}; pg. 280). However in the family of compact manifolds with NNSB, and even those with smooth boundary, and with $|Ric|\leq \Lambda_{0}$ and $diam_{g}(M)\leq D_{0}$ one cannot guarantee the existence
of ${\mathfrak N}_{0}:\mathbb{R}^{+}\rightarrow \mathbb{R}^{+}$ nor the precompactness of such family. Consider for instance the following example. For any $n\geq 2$ let $V_{n}=[1/n,1]\times \mathbb{S}^{1}$ be endowed with the flat metric $dx^{2}+n^{2}x^{2}d\varphi^{2}$ where $\varphi$ is the coordinate in the $\mathbb{S}^{1}$ factor (and recall that $\mathbb{S}^{1}$ has total length $2\pi$). For any $n$ the diameter of $V_{n}$ is less or equal than $2\pi+2$. On $M_{n}=V_{n}\times \mathbb{S}^{1}$ consider the flat product metric $g_{n}=dx^{2}+n^{2}x^{2}d\varphi^{2}+(1/n)^{2}d\theta^{2}$ where $\theta$ is the coordinate in the $\mathbb{S}^{1}$ factor defining $M_{n}$. Also, for any $n$, $diam_{g_{n}}(M_{n},g_{n})\leq 2\pi+2+1/n<10$. Despite of this and despite that the manifolds $(M_{n},g_{n})$ are flat, they do not collapse to a compact metric space (as $n\rightarrow \infty$). Even more we have that for any $1/2>\epsilon>0$ no pointed sequence $(\Omega_{n},g_{n},p_{n})$ of compact connected regions of $M_{n}$ with smooth boundary $\partial \Omega_{n}$, $\partial \Omega_{n}\subset {\mathcal T}_{d^{M_{n}}_{g_{n}}}(\partial M_{n},\epsilon)$, collapses to a compact metric space. This occurs even when $d^{M_{n}}_{g_{n}}(p_{n},\partial \Omega_{n})\geq 1/4$ (for instance).    

But any family ${\mathcal M}({\mathfrak N}_{0})$ satisfies the following kind of precompactness.
\begin{Proposition}\label{PCOM} Let $(M_{i},g_{i})$ be a sequence in a family ${\mathcal M}({\mathfrak N}_{0})$ and on every connected component of $M^{\circ}_{i}$ let $d_{i}=d^{M_{i}}_{g_{i}}$.  Then for every $1>\epsilon_{0}>0$ we have,
\begin{enumerate}
\item There are at most ${\mathfrak N}_{0}(\epsilon_{0},\epsilon_{0}/3)$ connected components $\breve{M}_{i}^{\circ}$ of $M_{i}^{\circ}$ intersecting 
$M_{i}\setminus {\mathcal T}_{d_{i}}(\partial M_{i},\epsilon_{0})$. 
\item For every sequence $\breve{M}_{i}^{\circ}$ of connected components of $M^{\circ}_{i}$ intersecting $M_{i}\setminus {\mathcal T}_{d_{i}}(\partial M_{i},\epsilon_{0})$, there is a subsequence (index again by ``i") such that $(\breve{M}_{i}^{\circ}\setminus {\mathcal T}_{d_{i}}(\partial M_{i},\epsilon_{0}),d_{i})$ converges in the GH-topology to a compact metric space $(X,d)$.
\end{enumerate}

\end{Proposition}  
\begin{Remark} 
Note that distances in $\breve{M}_{i}^{\circ}\setminus {\mathcal T}_{d_{i}}(\partial M_{i},\epsilon)$, which can be a connected set or not, are with respect to $d_{i}=d^{M_{i}}_{g_{i}}$.
\end{Remark}
{\bf Proof:} {\it Item 1}. By definition ${\mathfrak N}_{0}(\epsilon_{0},\epsilon_{0}/3)$ bounds from above the minimum number of balls of radius $\epsilon_{0}/3$ covering $M_{i}\setminus {\mathcal T}_{d_{i}}(\partial M_{i},\epsilon_{0})$. But given one such cover there must be at least one ball for every connected component $\breve{M}_{i}^{\circ}$ intersecting $M_{i}\setminus {\mathcal T}_{d_{i}}(\partial M_{i},\epsilon_{0})$. 
{\it Item 2}. From Definition \ref{DEF2} the function ${\mathfrak N}_{0}(\epsilon_{0},\epsilon_{1})$ as a function of $\epsilon_{1}$ and with $\epsilon_{0}$ fixed as in the hypothesis, bounds from above the minimum number of $d_{i}$-balls of radius $\epsilon_{1}$ covering $\breve{M}^{\circ}_{i}\setminus {\mathcal T}_{d_{i}}(\partial M_{i},\epsilon_{0})$. The Proposition follows from Lemma 1.9 in \cite{MR2243772} (pg. 280).\ep

\vs
\n In the example below we describe a nontrivial family of manifolds with boundary which are of great interest to us and lie in a class ${\mathcal M}({\mathfrak N}_{0})$. 
 
\vs
\n {\bf Example I.} Let $g$ be a (complete) Riemannian metric in $\mathbb{R}^{3}$. Suppose that $Ric_{g}\geq 0$ outside $B_{g}(o,r_{0})$ and that  $|Ric_{g}|\leq \Lambda_{0}/r^{2}$. Fix $0<c_{0}<c_{1}$. We claim that there is ${\mathfrak N}_{0}$ such that for any $\bar{r}$ with $\bar{r}\geq r_{1}(\Lambda_{0},r_{0},c_{0},c_{1})$, the Riemannian annuli (with NNSB) $(M_{\bar{r}},g_{\bar{r}})$,
\ben
M_{\bar{r}}:=\overline{A_{g}(c_{0}\bar{r},c_{1}\bar{r})}=\overline{B_{g_{\bar{r}}}(o,c_{1})\setminus \overline{B_{g_{\bar{r}}}(o,c_{0})}},\ \ g_{\bar{r}}:=\frac{1}{\bar{r}^{2}}g 
\een
lies in ${\mathcal M}({\mathfrak N}_{0})$. We show this in the following. From the {\it Ball-covering property} (\!\!\cite{MR1068127}, concretely Remark 2, pg. 215 (\footnote{In the Remark take $S=\overline{A_{g}(c_{0}\bar{r},c_{1}\bar{r})}$, $C_{0}=c_{1}/c_{0}$ and $\mu=c_{0}/3$.})) we know that for any $0<c_{0}<c_{1}$ there is $r_{1}$ and a number $n_{0}$ depending only on $c_{0}, c_{1}, r_{0}$ and $\Lambda_{0}$ such that for any $\bar{r}>r_{1}$, $n_{0}$ bounds from above the minimum number of $g_{\bar{r}}$-balls (in $\rt$) of radius $c_{0}/3$ covering the annulus  $M_{\bar{r}}$. Now, for any $g_{\bar{r}}$-ball with center in $M_{\bar{r}}$ and of $g_{\bar{r}}$-radius $c_{0}/3$, the minimum number of balls (in $\rt$) of $g_{\bar{r}}$-radii $\epsilon_{1}<c_{0}/4$ covering it (therefore having $Ric\geq 0$) is, by a simple application of the Bishop-Gromov volume comparison, bounded above by $(\frac{2c_{0}}{3})^{3}/(\epsilon_{1}^{3})$.  It follows that for any $1>\epsilon_{0}>2\epsilon_{1}$ with $\epsilon_{1}<c_{0}/4$, the minimum number of $g_{\bar{r}}$-balls (in $\rt$) of $g_{\bar{r}}$-radii $\epsilon_{1}$ covering $M_{\bar{r}}\setminus {\mathcal T}_{d_{\bar{r}}}(\partial M_{\bar{r}},\epsilon_{0})$ (here $d_{\bar{r}}=d^{M_{\bar{r}}}_{g_{r}}$) is bounded from above by $n_{0}c_{0}^{3}/\epsilon_{1}^{3}$. Therefore (recall Remark \ref{R1}) for any $\bar{r}>r_{1}$, $(M_{\bar{r}},g_{\bar{r}})$ belongs to ${\mathcal M}({\mathfrak N}_{0})$ where (when $1>\epsilon_{0}>2\epsilon_{1}>0$) ${\mathfrak N}_{0}(\epsilon_{0},\epsilon_{1})$ is defined as ${\mathfrak N}_{0}(\epsilon_{0},\epsilon_{1})=n_{0}c_{0}^{3}/\epsilon_{1}^{3}$ if $\epsilon_{1}<\min\{1/2,c_{0}/4\}$ and as ${\mathfrak N}_{0}(\epsilon_{0},\epsilon_{1})=4^{3}n_{0}$ if $\epsilon_{1}\in [\min\{1/2,c_{0}/4\},1/2)$. \et

\vs
\n {\bf Example II.} For any $D_{0}$, $\Lambda_{0}$ and $\delta_{0}$, there is $\mathfrak{N}_{0}(D_{0},\Lambda_{0},\delta_{0})$ such that for any $(M,g)$ Riemannian-manifold, with $|Ric_{g}|\leq \Lambda_{0}$, and connected compact region $\Omega\subset M$ with smooth boundary having
\ben
diam_{g}(\Omega)\leq D_{0},\text{ and } d^{M}_{g}(\partial \Omega,\partial M)> \delta_{0}, 
\een
the connected manifold (with NNSB) $\overline{{\mathcal T}_{d^{M}_{g}}(\Omega, \delta_{0})}$ lies is ${\mathcal M}(\mathfrak N_{0},\Lambda_{0})$. The proof is not difficult and is left to the reader.\et

\vs
It is instructive to go back and recall the discussion before the Proposition \ref{PCOM}. In there we presented a sequence $(M_{n},g_{n})$ which, in the light of Proposition \ref{PCOM}, did not belong to a single ${\mathcal M}(\mathfrak N_{0})$. Now, in the light of Example II, the manifolds $(M_{n},g_{n})$ (for all $n$) cannot be extended beyond their boundary to manifolds $(\bar{M}_{n},g_{n})$ with $|Ric_{g_{n}}|\leq \Lambda_{1}$ and $d^{\bar{M}_{n}}_{g_{n}}(M_{n},\partial \bar{M}_{n})\geq\delta_{0}>0$.  

As a consequence of Example II we have,

\vs
\n {\bf Example III}. Let $R_{0}>0$ and $\Lambda_{0}>0$ be given. Then, there is ${\mathfrak N}_{0}(R_{0},\Lambda_{0})$ 
such that any (closure of a) geodesic ball of radius $r_{0}\leq R_{0}$ inside a manifold $(M,g)$ with $|Ric|\leq \Lambda_{0}$, lies in ${\mathcal M}({\mathfrak N}_{0})$. To see this note that $B_{g}(p,r_{0})={\mathcal T}_{d^{M}_{g}}(B_{g}(p,r_{0}/2),r_{0}/2)$ and then use Example II. \et

\vs
We will denote by ${\mathcal M}({\mathfrak N}_{0},\Lambda_{0})$ the set of Riemannian three-manifolds (with NNSB) in the class ${\mathcal M}({\mathfrak N}_{0})$ and with $|Ric|\leq \Lambda_{0}$. In the Example I, the manifolds $(M_{\bar{r}},g_{\bar{r}})$ lie in  ${\mathcal M}({\mathfrak N}_{0},c_{0}^{-2}\Lambda_{0})$ where ${\mathfrak N}_{0}, c_{0}$ and $\Lambda_{0}$ are as in the example.
\begin{Definition}\label{ECCD} Let $(M,g)$ be a compact manifold (with NNSB). Let $0<\underline{\epsilon}<\overline{\epsilon}<1$. Then, a compact connected region 
$\Omega$ (with NNSB) is said to be an $(\underline{\epsilon},\overline{\epsilon})$-connected component of $M$ if $\partial \Omega\subset {\mathcal T}_{d_{g}^{M}}(\partial M,\overline{\epsilon})\setminus \overline{{\mathcal T}_{d_{g}^{M}}(\partial M,\underline{\epsilon})}$. The set of $(\underline{\epsilon},\overline{\epsilon})$-components of three-manifolds in a class ${\mathcal M}({\mathfrak N}_{0},\Lambda_{0})$ will be denoted by ${\mathcal M}^{\overline{\epsilon}}_{\underline{\epsilon}}({\mathfrak N}_{0},\Lambda_{0})$.
\end{Definition}
Thus when we write $(\Omega,g)\in {\mathcal M}^{\overline{\epsilon}}_{\underline{\epsilon}}({\mathfrak N}_{0},\Lambda_{0})$ we imply that $(\Omega,g)$ is the $(\overline{\epsilon},\underline{\epsilon})$-connected component of a $(M,g)\in {\mathcal M}({\mathfrak N}_{0},\Lambda_{0})$.

A sequence $(M_{i},g_{i})$ is {\it volume collapsing} if $Vol_{g_{i}}(M_{i})\rightarrow 0$. 
The following important Lemma is essentially Proposition 1.5 in \cite{MR1888088} (up to some modifications\footnote{\label{FONP}Unfortunately Proposition 1.5 in \cite{MR1888088} is stated without proof. An argumentative proof can be found in page 983 in \cite{MR1806984} (for the Lemma 1.4 in pg. 982 which is the equivalent to Proposition 1.5 in \cite{MR1888088}) but we were not able to check every claim in there, specially concerning the existence of $U_{i}$ (in the terminology of \cite{MR1806984}) with $\epsilon/2\leq dist(\partial U_{i},\partial \Omega_{i})\leq \epsilon$. The problems have to do with the fact that a priori the sequence $(D_{i},g_{i},x_{i})$ (in the terminology of \cite{MR1888088}) do not belong to any family ${\mathcal M}({\mathfrak N}_{0})$ and this may cause some inconveniences as indicated in the discussion before the Proposition \ref{PCOM}. It is essentially to avoid these inconveniences that we included the hypothesis that the sequence $(M_{i},g_{i})$ belongs a priori to some fixed family ${\mathcal M}({\mathfrak N}_{0})$. We would like to thank Michael Anderson for conversations on the Propositions 1.4 and 1.5 in \cite{MR1888088}.}) and with some additional information from \cite{MR950552}.
\begin{Lemma}\label{LL} Let $(M_{i},g_{i})$ be a volume-collapsing sequence in a ${\mathcal M}({\mathfrak N}_{0},\Lambda_{0})$ and such that for some $p_{i}\in M_{i}$ we have $d^{M_{i}}_{g_{i}}(p_{i},\partial M_{i})\geq \Gamma_{0}>0$. Then, for every $0<\underline{\epsilon}<\overline{\epsilon}<\min\{1,\Gamma_{0}/2\}$ there is a sequence $(\Omega_{i},g_{i})$ of $(\underline{\epsilon},\overline{\epsilon})$-connected components of $M_{i}$, with $p_{i}\in \Omega_{i}$, and a subsequence of it (indexed again by ``$i$") converging in the GH-topology to a space $(X,d)$ of one of the following two forms:
\begin{enumerate}
\item[\bf D1.] An interval $([0,\bar{x}],|\ |)$, with $Sing(X)=\emptyset$ or $Sing(X)=\{\bar{x}\}$, or,
\item[\bf D2.] A $C^{1,\beta}$-two-orbifold, with either $Sing(X)=\emptyset$ or $Sing(X)=\{\bar{x}_{1},\ldots,\bar{x}_{n}\}\subset X^{\circ}$.     
\end{enumerate}
Moreover (for $i\geq i_{0}$) {\bf I} and {\bf II} below hold.
\begin{enumerate}
\item[\bf I.] There are fibrations $f_{i}:\Omega_{i}\rightarrow X$, with asymptotically collapsing fibers $f^{-1}_{i}(x)$, such that,

- For {\bf D1}: $f_{i}:\Omega_{i}\setminus f^{-1}_{i}(Sing(X))\rightarrow X\setminus Sing(X)$ is a $\torus$-fibre-bundle and if $Sing(X)\neq \emptyset$ then $f_{i}^{-1}(\bar{x})\sim \torus/(\mathbb{S}^{1}\times {\mathbb{Z}_{q}})$, where the quotient is by a free action. 

- For {\bf D2}: $f_{i}:\Omega_{i}\setminus f^{-1}_{i}(Sing(X))\rightarrow X\setminus Sing(X)$ is a $\mathbb{S}^{1}$-fibre-bundle and if $Sing(X)\neq\emptyset$ then $f_{i}^{-1}(\bar{x}_{j})\sim \mathbb{S}^{1}/{\mathbb{Z}_{q_{j}}}$, where the quotient is by a free action. 

\item[\bf II.] There are finite coverings $\pi_{i}:\tilde{\Omega}_{i}\rightarrow \Omega_{i}$, such that

- For {\bf D1}: $(\tilde{\Omega}_{i},\tilde{g}_{i})$ converges in $C^{1,\beta}$ to a $\torus$-symmetric Riemannian manifold.

{\bf D2}: $(\tilde{\Omega}_{i},\tilde{g}_{i})$ converges in $C^{1,\beta}$ to a $\mathbb{S}^{1}$-symmetric Riemannian manifold.

In either case, for any $x\in X\setminus Sing(X)$, $\pi_{i}^{-1}(f^{-1}_{i}(x))$ converges in $C^{1}$ to the $\torus$ or $\mathbb{S}^{1}$ orbits.  

\end{enumerate}

\end{Lemma}

\vs
\n The fibrations $f_{i}$ have one more property \cite{MR873459}: for any neighborhood $W$ of $Sing(X)$ the map $f_{i}:f_{i}^{-1}(X\setminus W)\rightarrow X\setminus W$ is an almost Riemannian submersion, more precisely we have
\ben
e^{-o(i)}\leq |f_{i*}(V)|\leq e^{o(i)},\text{ where } o(i)\xrightarrow{i\rightarrow \infty} 0,
\een
and for any unit-norm $V$ perpendicular to the fibers.

\begin{Remark} We remark that the space $(\Omega_{i},g_{i})$ represents $(\Omega_{i},d^{\Omega_{i}}_{g_{i}})$ (see Sec. \ref{BNS}) rather than $(\Omega_{i},d^{M_{i}}_{g_{i}})$. Compare this with {\it item 2} in Proposition \ref{PCOM}. 
\end{Remark} 

Once one assumes that the sequence $(M_{i},g_{i})$ is in ${\mathcal M}({\mathfrak N}_{0},\Lambda_{0})$ the proof of Lemma \ref{LL} reduces to pointing to the appropriate reference in Fukaya's work. Here we overview why this is so. The proof itself is postponed to the Appendix. 

We introduce first a terminology. We say that two metric spaces $(Y,d_{Y})$ and $(Z,d_{Z})$ are {\it locally isometric under a homeomorphism} $\phi:Y\rightarrow Z$ if for all $y\in Y$ and $\phi(y)=z$ there are $\delta(y)$ and $\delta(z)$ such that $\phi:(B_{d_{Y}}(y,\delta(y)),d_{Y}) \rightarrow (B_{d_{Z}}(z,\delta(z)),d_{Z})$ is an isometry. Of course there are non-isometric metric spaces which are locally isometric \footnote{For instance compare the set $\{\varphi\in \mathbb{S}^{1},0< \varphi< 3\pi/4\}$ with the restriction of the standard metric in $\mathbb{S}^{1}$ and $((0,3\pi/4),|\ |)$.}. As a matter of fact if $(\Omega,g)\subset (M,g)$ then $(\Omega^{\circ},d^{\Omega}_{g})$ is locally isometric under the identity homeomorphism to $(\Omega^{\circ},d^{M}_{g})$, but they are not globally isometric in general. 

Suppose now that a sequence of compact boundaryless manifolds $(M_{i},g_{i})$ with uniformly bounded curvature and diameter
collapses to a metric space $(X,d)$ and suppose that $p_{i}\rightarrow x$. Let $exp:T_{p_{i}}M_{i}\rightarrow M_{i}$ be the exponential map and let $g_{i}(p_{i})$ be the metric $g_{i}$ on $T_{p_{i}}M_{i}$. Finally let $BT_{g_{i}(p_{i})}(p_{i},R_{0})$ be the $g_{i}(p_{i})$-ball of radius $R_{0}$ in $T_{p_{i}}M_{i}$. There is $R_{0}(\Lambda_{0})$ small enough, for which the map
\ben
exp: BT_{g_{i}(p_{i})}(p_{i},R_{0})\rightarrow B_{g_{i}}(p_{i},R_{0})
\een
is of maximal rank. Let $g_{i}^{*}$ be the pull-back metric. 
Then, Fukaya's technique to describe the space around $x$ (\!\!\cite{MR950552} Ch. 3), consists in working with $(BT_{g^{*}_{i}}(p_{i},R_{1}),g^{*}_{i})$, with $R_{1}\leq R_{0}$ small enough, and making the following observations\footnote{We do not comment here about some technical issues on smoothing.}
\begin{enumerate}
\item One can find a subsequence of it converging to a Riemannian manifold $(BT,g^{*})$ (\!\!\cite{MR950552}, pg. 9). 
\item For every $i$, $(B_{g_{i}}(p_{i},R_{1}/2),g_{i})$ is isometric to the quotient of the space $(BT_{g^{*}_{i}}(p_{i},R_{1}/2),g^{*}_{i})$ by an appropriate {\it local group}\footnote{See \cite{MR950552} and ref. therein.} of isometries $G_{i}$ and that $G_{i}$ converges to a local group $G$ (\!\! \cite{MR950552}, pg. 9) which is locally isomorphic to a Lie group (\!\!\cite{MR950552}, Lemma 3.1 in pg. 10). 
\item $(B_{d}(x,R_{1}/2),d)$ is locally isometric to $(BT(R_{1}/2),g^{*})/G$, where $BT(R_{1}/2)$ is the $g^{*}$-ball of radius $R_{1}/2$ in $BT$ (i.e the limit of $BT_{g^{*}_{i}}(p_{i},R_{1}/2)$) (\footnote{This is easy to check and is left to the reader.}).
\end{enumerate}
Thus by {\it item 3} to study locally the space $(X,d)$ around $x$ it is enough to study the limit spaces  $(BT(R_{1}/2),g^{*})/G$ and this is what is done in \cite{MR950552}. What is important to us about this conclusion is that one can study the collapse of manifolds with boundary as long as one works on a finite number of balls at a definite distance away from the boundary.   
This is essentially what is done in the proof of Lemma \ref{LL} in the Appendix and where the condition $(M_{i},g_{i})\in {\mathcal M}({\mathfrak N}_{0},\Lambda_{0})$ is used. 

\vs
We describe now a relevant application of Lemma \ref{LL} which will be of use to us in Proposition \ref{CRUC}. We describe it first in rough terms and then in a precise statement. Consider any solid torus with curvature bounded above by $\Lambda_{0}$ (fixed) and which is metrically close to an interval $I$ of length between $\infty\geq L_{0}>|I|\geq 1>0$ (with $L_{0}$ fixed) and with boundary metrically close to a point. Then, any curve ${\mathscr C}$ in its boundary, which is not a contractible to a point (from now on simply ``contractible") as a curve in the boundary, but that is contractible as a curve in the solid torus, must have length greater or equal than some $l_{0}(\Lambda_{0},L_{0})>0$. A proof of this phenomenon can be given along the following lines. Suppose that a curve ${\mathscr C}$ in the boundary of the solid torus $\Omega$, that is not a contractible curve as a curve in $\partial \Omega$ but is contractible as a curve in $\Omega$ has very small length. Then one can ``unwrap" $\Omega$, namely take a non-collapsed cover $\tilde{\Omega}$, which is also a solid torus. In particular $\partial \Omega$ is covered by a non-collapsed two-torus $\partial \tilde{\Omega}$.  But then the closed curve $\curve$, which is contractible in $\Omega$,  lifts to a closed, equal length and non-contractible curve $\tilde{\curve}$ in $\partial \tilde{\Omega}$. But there are no non-contractible curves $\tilde{\curve}$ in $\partial \tilde{\Omega}$ of very small length. This idea is made rigorous in the proof of Proposition \ref{PC2}.
This behavior is explicit in Example {\bf I.b} as we explain in what follows. In there the Riemannian solid tori $(M_{n},g_{n})$ are collapsing to a segment of length one. No matter the value of $n$, consider the $\curve_{0}=\pi_{n}(\tilde{\curve}_{0})$ where $\tilde{\curve}_{0}=\mathbb{S}^{1}\times \{1\}\subset \mathbb{B}^{2}\times \mathbb{S}^{1}$. The $g_{n}$-length of $\curve_{0}$ is equal to the length of $\tilde{\curve}_{0}$ and therefore equal to $2\pi$. Moreover any curve $\curve$ in $\partial M_{n}$ which is non-contractible as a curve in $\partial M_{n}$ but that is contractible as a curve in $M_{n}$ has length greater than that of $\curve_{0}$, i.e. $2\pi$. In other words, no matter the value of $n$, there are no such curves having a small length. 

We give an statement of what we described above in Proposition \ref{PC2}.
The statement is a bit more general than what was explained before as we do not make hypothesis on the boundary of the solid tori. For this reason too it is more general than what we will need in this article but it can be useful in other investigations. The proof is given in all detail partly to exemplify how the techniques apply.
\begin{Proposition}\label{PC2} For any $\Lambda_{0}$, $\delta_{0}<1/2$ and $L_{0}$ there is $\ell_{0}>0$ such that for any sequence $(\Omega_{i},g_{i})$ of solid tori inside a volume collapsing sequence of Riemannian manifolds $(M_{i},g_{i})$ with $|Ric_{g_{i}}|\leq \Lambda_{0}$, having
\begin{enumerate}
\item[\rm\bf Q0.] $rad_{g_{i}}(\Omega_{i})\geq 1$, $d^{M_{i}}_{g_{i}}(\partial \Omega_{i},\partial M_{i})\geq \delta_{0}>0$, and which is
\item[\rm\bf Q1.] Metrically collapsing to an interval $(I,|\ |)$, and,
\item[\rm\bf Q2.] Posses a sequence of closed curves ${\mathscr C}_{i}\subset \partial \Omega_{i}$ non-contractible in $\partial \Omega_{i}$ but contractible in $\Omega_{i}$ with $length_{g_{i}}({\mathscr C}_{i})\leq \ell_{0}$,
\end{enumerate}
we have $|I|\geq L_{0}$.
\end{Proposition}  
\begin{Remark} The hypothesis that the sequence $(M_{i},g_{i})$ is volume collapsing can be seen to be unnecessary.
\end{Remark}

\n {\bf Proof of Proposition \ref{PC2}:} For the proof it is worth to keep reference to the Figure \ref{FLL}. We will argue by contradiction. Suppose that there is $\Lambda_{0}$, $\delta_{0}<1/2$ and $L_{0}$ such that for every $m=1,2,3,\ldots$ there are sequences (in ``$i$") $(\Omega_{m,i},g_{m,i})\subset (M_{m,i},g_{m,i})$, where for every $m$, $(M_{m,i},g_{m,i})$ is a volume collapsing sequence of Riemannian manifolds with $|Ric_{g_{m,i}}|\leq \Lambda_{0}$, such that 
\begin{enumerate}
\item[$\rm\bf\bar{Q}0$.] $rad_{g_{m,i}}(\Omega_{m,i})\geq 1$, $d^{M_{m,i}}_{g_{m,i}}(\partial \Omega_{m,i},\partial M_{m,i})\geq \delta_{0}>0$, and which is,
\item[$\rm\bf\bar{Q}1$.] Metrically collapsing to an interval $(I_{m},|\ |)$, with $L_{0}\geq |I_{m}|$, and which 
\item[$\rm\bf\bar{Q}2$.] Posses a sequence of closed curves $\curve_{m,i}\subset \partial \Omega_{m,i}$ non-contractible in $\partial \Omega_{m,i}$ but contractible in $\Omega_{m,i}$ and of $length_{g_{m,i}}(\curve_{m,i})\leq 1/m$.
\end{enumerate}
Using that every sequence (in ``$i$") $(M_{m,i},g_{m,i})$ is volume collapsing and using ${\rm\bf\bar{Q}1}$, one can select for every $m$ an $i(m)$ such that 
\ben
Vol_{g_{m,i(m)}}(M_{m,i(m)})\leq 1/m,\text{ and  }dist_{GH}(\Omega_{m,i(m)}, I_{m})\leq 1/m. 
\een
In particular the sequence (in ``$m$") $(M_{m,i(m)},g_{m,i(m)})$ is volume collapsing. 
Also, because $rad_{g_{m,i(m)}}(\Omega_{m,i(m)})\geq 1$ and because of ${\rm\bf\bar{Q}1}$ we have $L_{0}\geq |I_{m}|\geq 1$ for every $m$ (\footnote{In general, if $(X_{m},d_{X_{m}})\xrightarrow{GH} (X,d_{X})$ and $d_{X_{m}}(x_{m},x'_{m})\geq \Gamma$ for all $m$, then there are $x$ and $x'$ in $X$ with $d_{X}(x,x')\geq \Gamma$ (use the definition of GH-convergence). On the other hand if $rad_{g_{m,i(m)}}(\Omega_{m,i(m)})\geq 1$ then there are $x_{m}$ and $x'_{m}$ in $\Omega_{m,i(m)}$ such that $d^{\Omega_{m,i(m)}}_{g_{m,i(m)}}(x_{m},x'_{m})\geq 1$.}). 
Therefore there is a subsequence of $(\Omega_{m,i(m)},g_{m,i(m)})$ (indexed again by ``$m$") metrically collapsing to an interval $I'$ with $L_{0}\geq |I'|\geq 1$. We continue working with this subsequence in what follows. 
This implies in particular that $diam_{g_{m,i(m)}}(\Omega_{m,i(m)})\leq D_{0}$ for some $D_{0}$ and for all $m$ (\footnote{In general if $(X_{m},d_{m})\xrightarrow{GH} (X,d)$ then there is $D_{0}$ such that $diam_{d_{X_{m}}}(X_{m})\leq D_{0}$ for all $m$ (use the definition of GH-convergence).}).  
Finally note, to be used below, that from ${\rm\bf\bar{Q}2}$ there is, for every $m$, a curve $\curve_{m,i(m)}\subset \partial \Omega_{m,i(m)}$ non-contractible in $\partial \Omega_{m,i(m)}$ but contractible in $\Omega_{m,i(m)}$ and of $length_{g_{i(m)}}(\curve_{m,i(m)})\leq 1/m$.

By the Example II, if we let $M'_{m}={\mathcal T}_{d_{m}}(\Omega_{m,i(m)}, \delta_{0})$ with $d_{m}=d_{g_{m,i(i)}}^{M_{m,i(m)}}$, then $(M'_{m},g_{m,i(m)})$ lies in ${\mathcal M}({\mathfrak N}_{0},\Lambda_{0})$ for some ${\mathfrak N}_{0}(D_{0},\Lambda_{0},\delta_{0})$. On the other hand as $M'_{m}\subset M_{m,i(m)}$ then $(M'_{m},g_{m,i(m)})$ is also a volume collapsing sequence.
Hence, by Lemma \ref{LL}, one can find a sequence of $(\delta_{0}/4,\delta_{0}/2)$-connected components of  $M'_{m}$ containing $\Omega_{m,i(m)}$, to be denoted by $\hat{\Omega}_{m}$, and having a subsequence (indexed again by ``$m$") metrically collapsing to an interval $\hat{I}$ containing $I'$. We continue using this subsequence in what follows.
For the sake of concreteness assume that $\hat{I}$ is the interval $[0,|\hat{I}|]$.

Consider the fibrations $f_{m}:\hat{\Omega}_{m}\rightarrow \hat{I}$ as is explained in Lemma \ref{LL}. As $m\rightarrow \infty$, the fibers $f^{-1}_{m}(x)$ collapse to a point and so does $\partial \hat{\Omega}_{m}=f^{-1}_{m}(0)$ to the point $0$ in $\hat{I}$. Observe that the right point of $I'$ must be the right point of $\hat{I}$, that is $|\hat{I}|$, and therefore it is a singular point, namely $Sing(\hat{I})=\{|\hat{I}|\}$. 
We observe too that from the very definition of $M'_{m}$ we have, for every $q\in \partial \Omega_{m,i(m)}$, 
$d^{\hat{\Omega}_{m}}_{g_{m,i(m)}}(q,\partial \hat{\Omega}_{m})< \delta_{0}<1/2$ (\footnote{Note for this that for any $q\in \partial \Omega_{m,i(m)}$ we must have $B_{g_{m,i(m)}}(q,\delta_{0})\cap \partial \hat{\Omega}_{m}\neq \emptyset$, because $\hat{\Omega}_{m}$ is a $(\delta_{0}/4.\delta_{0}/2)$-c.c.}).
It follows from this that for $m\geq m_{0}$ with $m_{0}$ big enough 
(i) $\partial \Omega_{m,i(m)}\subset f_{m}^{-1}([0,1/2])$, (ii) $f_{m}^{-1}(1/2)$ lies in the interior of $\Omega_{m,i(m)}$ and (iii) $f^{-1}_{m}(0)$ lies in the exterior of $\Omega_{m,i(m)}$. In this way $\partial \Omega_{m,i(m)}$ separates $f^{-1}_{m}([0,1/2])$, which is diffeomorphic to $\torus\times [0,1/2]$, into two connected components. This implies\footnote{This is a simple exercise in topology (use Alexander's theorem in \cite{Alexander1} for two-tori in $\mathbb{S}^{3}$).} that $\partial \Omega_{m,i(m)}$ is isotopic to $f^{-1}_{m}(x)$ for any $x\in [0,1/2]$. In particular if $\curve_{m,i(m)}$ is non-contractible in $\partial \Omega_{m,i(m)}$ then it is also non-contractible in $f^{-1}_{m}([0,1/2])$. 
Moreover, by Lemma \ref{LL}, there is a subsequence (indexed again by ``$m$") and coverings $\pi_{m}:\tilde{\hat{\Omega}}_{m}\rightarrow \hat{\Omega}_{m}$ such that $(\tilde{\hat{\Omega}}_{m},\tilde{g}_{m,i(m)})$ converges in $C^{1,\beta}$ to a $\torus$-symmetric metric on $\mathbb{B}^{2}\times \mathbb{S}^{1}$ and $(\pi_{m}^{-1}(f^{-1}_{m}([0,1/2])),\tilde{g}_{m,i(m)})$ converges in $C^{1,\beta}$ to a $\torus$-symmetric metric on $\torus\times \mathbb{I}$.
For this reason there are $m_{1}$ and $\ell_{1}$, such that for any $m\geq m_{1}$ any non contractible closed curve in $(\pi_{m}^{-1}(f^{-1}_{m}([0,1/2])),\tilde{g}_{m,i(m)})$ has length greater or equal than $\ell_{1}$. 
But for every $m$, the curve $\curve_{m,i(m)}$ is closed and contractible in $\Omega_{m,i(m)}$ and thus contractible also in $\tilde{\hat{\Omega}}_{m}$. Therefore its lift
$\tilde{\curve}_{m,i(m)}$ to $\pi_{m}^{-1}(f_{m}^{-1}[0,1/2]))\subset\tilde{\hat{\Omega}}_{m}$ is also closed and has the same length, which, as was observed above, is less or equal than $1/m$. If $m\geq \max \{m_{1},2/\ell_{1}\}$ then $length_{g_{i(m)}}(\curve_{m,i(m)})\leq \ell_{1}/2$ which is not possible.\ep

\vs
\begin{figure}[h]
\centering
\includegraphics[width=12cm,height=7cm]{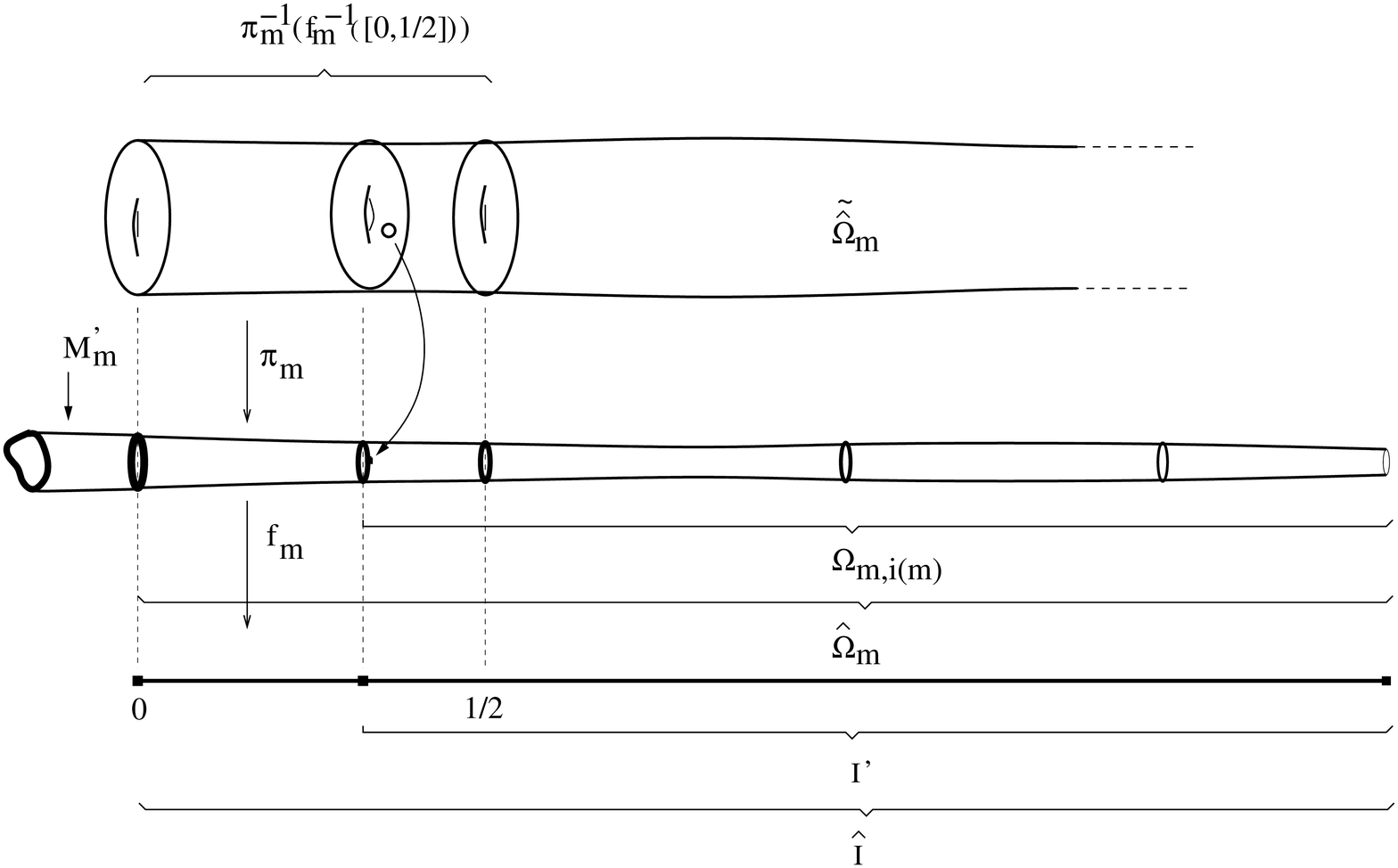}
\caption{A representation of the argument given in the proof of Proposition \ref{PC2}. The little curve in the cover manifold represents the lift $\tilde{\curve}_{m,i(m)}$ of $\curve_{m,i(m)}$. If $m\geq \max\{m_{1},2\ell_{1}\}$ the length of $\tilde{\curve}_{m,i(m)}$ would be too small to be non-contractible in $\pi_{m}^{-1}(f^{-1}_{m}([0,1/2]))$.}
\label{FLL}
\end{figure} 
\subsection{ A special annuli decomposition.}
\label{SAD}
The results of the previous section allow us to show the existence of an annuli decomposition with special properties.
\begin{Proposition}\label{SADP} Let $g$ be a complete metric in $\rt$ with
\ben
|Ric_{g}|\leq \frac{\Lambda_{0}}{r^{2}}, \text{ and }\lim_{\bar{r}\uparrow \infty} \frac{Vol_{g}(B_{g}(o,\bar{r}))}{\bar{r}^{3}}=0.
\een
Then, there is an annuli decomposition ${\mathcal U}$ with the following properties: for  every $\epsilon>0$ there is $k(\epsilon)$ such that for any $k\geq k(\epsilon)$ every piece $(U_{k,l},g_{k})$ is $\epsilon$-close in the GH-metric to a space $X_{k,l}$ of one of the following two forms, 
\begin{enumerate}
\item[$\rm\bf \tilde{D}1$.] An interval, in which case $U_{k,l}$ is either diffeomorphic to $\torus\times \mathbb{I}$ or a solid torus $\mathbb{B}^{2}\times \mathbb{S}^{1}$, or,
\item[$\rm\bf \tilde{D}2$.] A two-orbifold, in which case $U_{k,l}$ is diffeomorphic to a Seifert manifold with at least one boundary component. 
\end{enumerate}
There are fibrations $f_{k,l}:U_{k,l}\rightarrow X_{k,l}$, such that  for any $k\geq k(\epsilon)$ the fibers $f^{-1}_{k,l}(x)$, which are diffeomorphic either to $\torus$ or $\mathbb{S}^{1}$, are $\epsilon$-collapsed. Moreover
\begin{enumerate}
\item[$\rm\bf \tilde{I}1$.] In case {$\rm\bf \tilde{D}1$}, either $Sing(X_{k,l})$ is empty or is one of the extreme points of the interval. In addition, for any non singular point $x$, the fiber $f^{-1}_{k,l}(x)$ is diffeomorphic to $\torus$ and if $x$ is a singular point then $f^{-1}_{k,l}(x)$ is diffeomorphic to $\mathbb{S}^{1}$.   
\item[$\rm\bf \tilde{I}2$.] In case {$\rm\bf \tilde{D}2$}, the fibers $f^{-1}_{k,l}(x)$, which are all diffeomorphic to $\mathbb{S}^{1}$, are the fibers of the Seifert-fibration.
\end{enumerate}
%
%
%
\end{Proposition} 

\vs
\n Before going into the proof we introduce some notation. For every $k$ we define the scaled metric 
\be\label{NKG}
g_{k}=\frac{1}{10^{2k}}g
\ee
Therefore $A_{g}(10^{n_{1}+k},10^{n_{2}+k})=A_{g_{k}}(10^{n_{1}},10^{n_{2}})$ which to simplify notation we will write simply as
$A_{k}(10^{n_{1}},10^{n_{2}})$. 
We say that a set of embedded two manifolds $\{S_{k,j},j=1,\ldots,j(k)\}$ is {\it ``a cut of $\rt$ along the annulus $A_{k}(10^{-1},1)$}" if 
\begin{enumerate}
\item $S_{k,j}\subset A_{k}(10^{-1},1)$ for all $j=1,\ldots,j(k)$, and,
\item Every curve $\alpha:[0,1]\rightarrow \rt$ with $\alpha(0)\in B_{g_{k}}(o,10^{-1})$ and $\alpha(1)\in \big(\rt\setminus \overline{B_{g_{k}}(o,1)}\big)$ intersects at least one of the $S_{k,j}$'s, and,
\item The {\it item 2} does not hold if one deletes one of the $S_{k,j}$'s from the set. 
\end{enumerate}
Observe that if a set of manifolds $\{\bar{S}_{k,j}\}$ enjoy {\it item 1} and {\it item 2}, then one can remove, if necessary, some elements of the set to satisfy also {\it item 3}. Also, any surface $S_{k,j}$ of a ``cut" is necessarily the boundary of two connected components of $\rt\setminus \bigcup_{j=1}^{j=j(k)} S_{k,j}$, one intersecting $B_{g_{k}}(o,10^{-1})$ and the other intersecting $\rt\setminus \overline{B_{g_{k}}(o,1)}$ (\footnote{A connected component cannot intersect $B_{g_{k}}(o,10^{-1})$ and  $\rt\setminus \overline{B_{g_{k}}(o,1)}$ simultaneously, otherwise {\it item 2} is violated. Also, if a connected component does not intersect any of them then one can remove any boundary component still satisfying {\it item 2} and therefore violating {\it item 3}.}).  

\vs
\n {\bf Proof:} For the proof it may be worth to keep in mind the Figure \ref{UU}. As explained in the Example I in Section \ref{VCTM}, the spaces $(\overline{A_{k}(10^{-2},10^{4})},g_{k})$ lie in ${\mathcal M}({\mathfrak N}_{0},10^{2}\Lambda_{0})$ for some $k$-independent ${\mathfrak N}_{0}$. Moreover for any $p\in \overline{A_{k}(1,10)}$ we have $d^{\rt}_{g_{k}}(p,\partial \overline{A_{k}(10^{-2},10^{4})})>1/2$. Granted these two facts we can use then Lemma \ref{LL}  to obtain with no difficulty that: 
%
%
%
%
%

{\it There is a set $\{\bar{U}_{k,j},j=1,\ldots,j(k); k=k_{0},k_{0}+2,\ldots\}$ of (for each $k$) $(10^{-2}/2,10^{-2})$-connected components of $\overline{A_{k}(10^{-2},10^{4})}$ with the following properties. 
\begin{enumerate}
\item The set $\{\bar{U}_{k,j},j=1,\ldots,j(k)\}$ covers $\overline{A_{k}(10^{-1},10^{3})}$ for every $k=k_{0},k_{0}+2,\ldots$.
\item There are intervals or two-orbifolds, to be denoted by $\bar{X}_{k,l}$, and for every $m=1,2,3,\ldots$ there is $k_{m}$, such that if $k\geq k_{m}$ then $(\bar{U}_{k,j},g_{k})$ is $1/m$-close in the GH-metric to $\bar{X}_{k,l}$. 
\item There are fibrations $\bar{f}_{k,j}:\bar{U}_{k,j}\rightarrow \bar{X}_{k,j}$, with the properties ${\rm\bf \tilde{I}1}$ and ${\rm\bf \tilde{I}2}$,
such that if $k\geq k_{m}$ their fibers are $1/m$-collapsed. 
\end{enumerate}}    
\n The fibrations $\bar{f}_{k,l}:\bar{U}_{k,l}\rightarrow \bar{X}_{k,l}$ can be chosen in such a way that if $\bar{U}_{k,l}$ and $\bar{U}_{k',l'}$ overlap and have fibers of the same dimension, namely both have fibers of dimension one or both have fibers of dimension two, then the foliations of fibers coincide on the overlap, while if one has fibers of dimension one and the other of dimension two, then fibers of dimension one are included in fibers of dimension two. For this the reader can consult the geometric construction of the fibrations in \cite{MR873459} and \cite{MR950552}.

The desired sets $U_{k,l}$ of the annuli decomposition will be defined below simply as regions of the sets $\bar{U}_{k,j}$ appropriately ``cut out along the annuli $A_{k}(10^{-1},1)$" using the fibers of the fibrations $\bar{f}_{k,l}$. This has to be done in such a way to satisfy {\it items 1-5} of the definition of annuli decompositions. Once this is performed the
fibrations $f_{k,l}:U_{k,l}\rightarrow X_{k,l}$ are defined by $f_{k,l}:=\bar{f}_{k,l}|_{U_{k,l}}:U_{k,l}\rightarrow X_{k,l}:=\bar{f}_{k,l}(U_{k,l})$, where $\bar{U}_{k,j}$ is that piece containing $U_{k,l}$ and $\bar{f}_{k,l}:\bar{U}_{k,l}\rightarrow \bar{X}_{k,l}$ its fibration. We explain how the regions $U_{k,l}$ are constructed in what follows. 

Fix a value of $k$ in $\{k_{0},k_{0}+2,\ldots\}$. Then, on those $\bar{X}_{k,j}$ which are an interval select a set of points $\bar{x}_{j}$  and then on those $\bar{X}_{k,j}$ which are a two-orbifold select a set of (disjoint) closed curves denoted by $\bar{\curve}_{k,j,i}$, such that the set of tori
$\{S_{k,j}\}:=\{\bar{f}_{k,j}^{-1}(\bar{x}_{k,j}), \bar{f}^{-1}_{k,j,i}(\bar{\curve}_{k,j,i}), \text{ all } i,j,k\}$ is a ``cut of $\rt$ along the annulus $A_{k}(10^{-1},1)$" as defined before the start of the proof.  

Now, for every $k$ in $\{k_{0}, k_{0}+2,\ldots\}$ let $\hat{\mathcal U}_{k}$ be the set of compact connected regions of $\rt$ with boundary components in $\{S_{k,j}, S_{k+2,j'},\text{ all } j\text{ and } j'\}$. As mentioned before the start of the proof every $S_{k,l}$ is the boundary of two of such regions: one in $\hat{\mathcal U}_{k}$ and intersecting $\rt\setminus \overline{B_{g_{k}}(o,1)}$, denoted from now on by $\hat{U}^{+}(S_{k,j})$, and the other in $\hat{\mathcal U}_{k-2}$ and intersecting $B_{g_{k}}(o,10^{-1})$, denoted from now on by $\hat{U}^{-}(S_{k,j})$. Moreover we have the following two properties.
\begin{enumerate}
\item For every $S_{k,j}$, the piece $\hat{U}^{-}(S_{k,j})$ is equal to a piece $\hat{U}^{+}(S_{k-2,j'})$ for some $S_{k-2,j'}$, but not necessarily every piece $\hat{U}^{+}(S_{k,j})$ is a piece $\hat{U}^{-}(S_{k+2,j'})$ (\footnote{To see $1$ observe that if a piece $\hat{U}^{-}(S_{k,j})$ is not a $\hat{U}^{+}(S_{k-2,j'})$-piece then it does not have a boundary component in $\{S_{k-2,j},\text{ all } j\}$. Therefore, as $\hat{U}^{-}(S_{k,j})$ must intersect $B_{g_{k}}(o,10^{-1})$ there is a point $p\in \hat{U}^{-}(S_{k,j})$ with $d^{\rt}_{g_{k}}(p,o)<1$. Now, any length minimizing geodesic segment joining $p$ and $o$ must intersect $\partial \hat{U}^{-}(S_{k,j})$, say at $q$. But $\partial \hat{U}^{-}(S_{k,j})\subset A_{k}(1,10)$ and therefore $d^{\rt}_{g_{k}}(q,o)>1$ which is impossible as $q$ belongs to the length minimizing geodesic.}). 
\item Every piece $\hat{U}^{+}(S_{k,j})$ is included in $A_{k}(10^{-1},10^{2})$ (\footnote{To see this use directly {\it item 2} of the definition of ``cut".}).  
\end{enumerate}
\n We define now the pieces $U_{k_{0},l}$: 
\begin{itemize}
\item Redefine the set $\hat{\mathcal U}_{k_{0}}$ by eliminating from it those pieces $\hat{U}^{+}(S_{k_{0},j})$ which are not equal to a piece $\hat{U}^{-}(S_{k_{0}+2,j'})$. 
\item Every $\hat{U}^{+}(S_{k_{0}+2,j})$ which is not a $\hat{U}^{-}(S_{k_{0}+4,j'})$ is glued to those pieces in $\hat{\mathcal U}_{k_{0}}$ sharing a boundary component with it. We define the resulting manifold as one of the $U_{k_{0},l}$'s. 
The other pieces $U_{k_{0},l}$ are defined as those in $\hat{\mathcal U}_{k_{0}}$ which were not glued to a $\hat{U}^{+}(S_{k_{0}+2,j})$ as was explained. In either case every piece $U_{k_{0},l}$ is included in $A_{k_{0}}(10^{-1},10^{3})$ (use {\it 2}) and in one and only one of the $\bar{U}_{k_{0},j}$'s, which as was explained above is used to defined the fibrations $f_{k_{0},l}:U_{k_{0},l}\rightarrow X_{k_{0},l}$.  
\end{itemize}

\n We define next the pieces $U_{k_{0}+2,l}$:
\begin{itemize}
\item Redefine the set $\hat{\mathcal U}_{k_{0}+2}$ by eliminating from it those pieces $\hat{U}^{+}(S_{k_{0}+2,j})$ which are not equal to a piece $\hat{U}^{-}(S_{k_{0}+4,j'})$ and which, as was explained before, were glued to pieces in $\hat{\mathcal U}_{k_{0}}$ to form some of the pieces $U_{k_{0},l}$. 
\item Every $\hat{U}^{+}(S_{k_{0}+4,j})$ which is not a $\hat{U}^{-}(S_{k_{0}+6,j'})$ is glued to those pieces in $\hat{\mathcal U}_{k_{0}+2}$ sharing a boundary component with it. We define the resulting manifold as one of the pieces $U_{k_{0}+2,l}$. The other pieces $U_{k_{0}+2,l}$ are defined as those in $\hat{\mathcal U}_{k_{0}+2}$ which were not glued to a $\hat{U}^{+}(S_{k_{0}+4,j})$ as was explained.  In either case every piece $U_{k_{0}+2,l}$ is included in $A_{k_{0}+2}(10^{-1},10^{3})$ (use {\it 2}) and in one and only one of the $\bar{U}_{k_{0}+2,j}$'s, which as was explained above is used to defined the fibrations $f_{k_{0}+2,l}:U_{k_{0}+2,l}\rightarrow X_{k_{0}+2,l}$.  
\end{itemize}
To define the pieces $U_{k_{0}+4,l}$, $U_{k_{0}+6,l}$ and so on, proceed in the same way as the pieces $U_{k_{0}+2,l}$ were defined. It is straightforward to check that the family ${\mathcal U}=\{U_{k,l}\}$ thus defined  satisfies the Definition \ref{D1} of annuli decomposition.\ep 
\section{ Proof of Theorem \ref{MT}.}\label{PMT}
We will work in this section with the annuli decomposition defined in the previous section. We already defined in Section \ref{AD} the set ${\mathcal N}$ of boundary components of ${\mathcal U}$ which we will denote here generically by $T^{2}$ (instead of $S$ because they are tori). We also defined the subclass ${\mathcal N}^{o}$ as those tori $T^{2}$ in ${\mathcal N}$ for which $o\in M(T^{2})$ and observed that they were linearly ordered, i.e. ${\mathcal N}^{o}=\{T^{2o}_{0},T^{2o}_{1},\ldots\}$, with $T^{2o}_{i}\ll T^{2o}_{i'}$ if $i<i'$. For later convenience we further divide ${\mathcal N}\setminus {\mathcal N}^{o}$ into two subclasses denoted by ${\mathcal N}^{\blacklozenge}$ and ${\mathcal N}^{\lozenge}$; ${\mathcal N}^{\blacklozenge}$ (resp. ${\mathcal N}^{\lozenge}$) is defined as the set of tori in ${\mathcal N}\setminus {\mathcal N}^{\circ}$ for which $M(T^{2})$ is a solid torus (resp. not a solid torus).
Tori in ${\mathcal N}^{\blacklozenge}$ (resp. ${\mathcal N}^{\lozenge}$) will denoted as $T^{2\blacklozenge}$ (resp. $T^{2\lozenge}$). For every $T^{2}$ in ${\mathcal N}$ there is a unique piece $U_{k,l}$ (including the possibility of $U_{k_{0}-2}$) such that 
$T^{2}\in U_{k,l}$ and $U_{k,l}\subset M(T^{2})$. In this way the indexes $k,l$ are univocally defined and we can write $k(T^{2}), l(T^{2})$. We will continue using the notation (\ref{NKG}) in particular we will use $g_{k(T^{2}}$.

The following proposition is crucial for the proof of the Theorem \ref{MT}. Observe that the statement is suitable to be used in an iterative argument as will be the case when we use it in the proof of Theorem \ref{MT}. 
\begin{Proposition}\label{CRUC} There exits $\epsilon^{*},\ell^{*},k^{*}$ such that if for a $\st_{1}\in {\mathcal N}^{\blacklozenge}$ with $k(\st_{1})\geq k^{*}$ we have

\begin{enumerate}
\item[{\bf H1}.]$(U_{k(\st_{1}),l(\st_{1})},g_{k(\st_{1})})$ is $\epsilon^{*}$-close in the GH-metric to an interval, and,
\item[{\bf H2}.] There is a curve ${\mathscr C}_{1}\subset \st_{1}$ non-contractible in $\st_{1}$ but contractible in $M(\st_{1})$ such that $length_{g_{k(\st_{1})}}({\mathscr C}_{1})\leq \ell^{*}$,   
\end{enumerate}
then, $U_{k(\st_{1}),l(\st_{1})}$ is not the only $U_{k,l}$-piece of $M(\st_{1;m})$, and, if we denote by $\st_{2}$ the second boundary component of $U_{k(\st_{1}),l(\st_{1})}$, we have
\begin{enumerate}
\item[{\bf C1}.] $(U_{k(\st_{2}),l(\st_{2})},g_{k(\st_{2})})$ is $2\epsilon^{*}/3$-close in the GH-metric to an interval, and,
\item[{\bf C2}.] There is a curve ${\mathscr C}_{2}\subset \st_{2}$ non-contractible in $\st_{2}$ but contractible in $M(\st_{2})$ such that $length_{g_{k(\st_{2})}}({\mathscr C}_{2})\leq 2\ell^{*}/3$.
\end{enumerate}

\end{Proposition}

\vs
\n {\bf Proof:} By contradiction assume that for every $\epsilon^{*}_{m}=1/m, \ell^{*}_{m}=1/m$ and $k^{*}_{m}=m$, $m=1,2,3,\ldots$, there is $\st_{1;m}\in {\mathcal N}^{\blacklozenge}$ with $k(\st_{1;m})\geq k^{*}_{m}$ such that 
\begin{enumerate}
\item[$\rm\bf\bar{H}1$.] $(U_{k(\st_{1;m}),l(\st_{1;m})},g_{k(\st_{1;m})})$ is $\epsilon^{*}_{m}$-close in the GH-metric to an interval, and,
\item[$\rm \bf \bar{H}2$.] There is a curve ${\mathscr C}_{1;m}\subset \st_{1;m}$ non-contractible in $\st_{1;m}$ but contractible in $M(\st_{1;m})$ such that $length_{g_{k(\st_{1;m})}}({\mathscr C}_{1;m})\leq \ell^{*}_{m}$,   
\end{enumerate}
but that, if it is not that $U_{k(\st_{1;m}),l(\st_{1;m})}= M(\st_{1;m})$, then, after denoting by $\st_{2;m}$ the second boundary component of $U_{k(\st_{1;m}),l(\st_{1;m})}$, one of the following two assertions does not hold 
\begin{enumerate}
\item[${\rm \bf \bar{C}1}$.] $(U_{k(\st_{2;m}),l(\st_{2;m})},g_{k(\st_{2;m})})$ is $2\epsilon^{*}_{m}/3$-close in the GH-metric to an interval,
\item[$\rm \bf \bar{C}2$.] There is a curve ${\mathscr C}_{2;m}\subset \st_{2;m}$ non-contractible in $\st_{2;m}$ but contractible in $M(\st_{2;m})$ such that $length_{g_{k(\st_{2;m})}}({\mathscr C}_{2;m})\leq 2\ell^{*}_{m}/3$.
\end{enumerate}
We will show that this leads to an impossibility. Such impossibility will come directly as the result of proving the following three steps.

\vspace{0.4cm}
$\bullet$ {\bf Step A}. {\it Let $\st_{1;m}$ be a sequence satisfying ${\rm \bf \bar{H}1}$ and ${\rm \bf \bar{H}2}$. Then }
\ben
rad_{g_{k(\st_{1,m})}}(M(\st_{1;m}))\xrightarrow{m\rightarrow \infty} \infty.
\een

\n Step {\bf A} shows that there is $m_{1}$ such that for every $m\geq m_{1}$, $U_{k(\st_{1;m}),l(\st_{1;m})}$ is not the only piece of $M(\st_{1;m})$ (because if it is so then, by {\it item 1} of Definition \ref{D1}, $rad_{g_{k}(\st_{1;m})}(M(\st_{1;m}))\leq 10^{3}$). 
%
The statement of Step {\bf B} below assumes $m\geq m_{1}$.

\vs
$\bullet$ {\bf Step B}. {\it $(m\geq m_{1})$. Let $\st_{2;m}$ be the second component of $U_{k(\st_{1;m}),l(\st_{1;m})}$. Then there is a  covering sequence to a subsequence of 
\ben
(U_{k(\st_{1;m}),l(\st_{1;m})}\cup U_{k(\st_{2;m}),l(\st_{2;m})}, g_{k(\st_{1;m})})
\een
converging in $C^{1,\beta}$ to a flat $\torus$-symmetric metric product on $\mathbb{T}^{2}\times I_{1,2}$ for some interval $I_{1,2}$. That is, the limit metric on $\mathbb{T}^{2}\times I$ is of the form $dx^{2}+\tilde{h}_{0}$ with $\tilde{h}_{0}$ an ($x$-independent) $\torus$-symmetric metric on $\mathbb{T}^{2}$.}

\vs
$\bullet$ {\bf Step C}. {\it There is $m_{2}\geq m_{1}$ such that for all $m\geq m_{2}$ and $m$ in the subsequence of Step {\bf B}, ${\rm \bf \bar{C}1}$ and ${\rm \bf \bar{C}2}$ hold.}

\vspace{0.4cm}
\n From now until the end of the proof of the Proposition and to simplify notation we let
\begin{gather*}
U_{1;m}=U_{k(\st_{1;m}),l(\st_{1;m})},\ \ M_{1;m}=M(\st_{1;m}),\\ 
g_{1;m}=g_{k(\st_{1;m})},\ \ k_{1;m}=k(\st_{1;m})
\end{gather*}

{\bf\small Proof of Step A}. Assume on the contrary $rad_{g_{1;m}}(M_{1;m})\leq R_{0}$. Then, it is simple to see
\footnote{For any $p\in M_{1}(m)$, $d^{\rt}_{g_{1;m}}(p,o)$ is less or equal than
$d^{\rt}_{g_{1;m}}(p,\st_{1;m})+d^{\rt}_{g_{1;m}}(\st_{1;m})+d^{\rt}_{g_{1;m}}(\st_{1;m},o)$
which is less or equal than $R_{0}+diam_{g_{1;m}}(\st_{1;m})+1$. But the $g_{1;m}$-diameter of $\st_{1;m}$ tends to zero (by ${\bf\rm \bar{H}1}$) and so we can assume that it is less or equal than some $D_{0}$.}
that $M_{1;m}$ must be a subset of an annulus $A(10^{k_{1;m}-1},10^{k_{1;m}+k_{\bullet}})$ for some $k_{\bullet}>0$ independent of $m$. On the other hand $rad_{g_{1;m}}(M_{1;m})\geq 10^{2}-10^{-1}> 90$ (because of {\it item 2} of Definition \ref{D1} applied to $U_{1;m}$). Under these hypothesis we obtain
\begin{enumerate}
\item (using ${\rm\bf \bar{H}1}$) A subsequence of the sequence of solid tori $(M_{1;m},g_{1;m})$ (indexed still by ``$m$") metrically collapses to a compact interval\footnote{It must converge to an interval and not a two-orbifold (the only two options) because $(M_{1;m},g_{1;m})$ contains $(U_{1;m},g_{1;m})$ which by ${\rm\bf \bar{H}1}$ converges to an interval.} $I$ of length $|I|$ greater or equal than $90$, and,
\item (using ${\rm\bf \bar{H}2}$) For every $\ell_{0}$ we have, $\lim_{m\rightarrow \infty} length_{g_{1;m}}({\mathscr C}_{1;m})\leq \ell_{0}$. 
\end{enumerate}  
We can then apply Proposition \ref{PC2} (\footnote{To apply Proposition \ref{PC2} use as $M_{m}$ in its statement the manifold $M_{m}:={\mathcal T}_{d^{\rt}_{g_{1;m}}}(M_{1;m},10^{-2})$. It is direct that $(M_{m},g_{1;m})$ is a volume collapsing sequence.}) to conclude that $|I|\geq L_{0}$ for any $L_{0}$ and therefore that $|I|=\infty$, contradicting the compactness of the interval $I$. \et

\vs
We recount a little the setup and terminology before we go into Step {\bf B}. Let $\st_{2;m}$ be the second boundary component of $U_{1;m}$ and let $U_{2;m}:=U_{k(\st_{2;m}),l(\st_{2;m})}$ be the $U_{k,l}$-piece, other than $U_{1;m}$, having $\st_{2;m}$ as a boundary component. Of course $k_{2;m}:=k(\st_{2;m})=k(\st_{1;m})+2=k_{1;m}+2$. 
Following the same pattern of notation as before we let
\ben
U_{1,2;m}=U_{1;m}\cup U_{2;m},\ \
g_{2;m}=g_{k(\st_{2;m})}
\een

{\small\bf Proof of Step B}. 
%
%
To this end first note that $(U_{1,2;m},g_{1;m})$ collapses metrically to an interval\footnote{Again, this is so because $(U_{1;m},g_{1;m})$, with $U_{1}(m)\subset U_{1,2;m}$,  collapses metrically to an interval.} to be denoted by $I_{1,2}$; $(U_{1;m},g_{1;m})$ collapses to $I_{1}$ and $(U_{2;m},g_{2;m})$ collapses to $I_{2}$, and we have $I_{1,2}=I_{1}\cup I_{2}$ and $|I_{1,2}|=|I_{1}|\cup |I_{2}|$. Without loss of generality we assume that $\st_{1;m}$ collapses to the left boundary point of the interval $I_{1}$ (or, the same, of $I_{1,2}$) as an interval in $\mathbb{R}$.
Further, following Proposition \ref{SADP} and Lemma \ref{LL} (see also Proposition \ref{AP1} in the Appendix for a technical point on the explicit form of the limit), there is a subsequence (indexed again by $m$) and a covering sequence $\pi_{m}:\tilde{U}_{1,2;m}\rightarrow U_{1,2;m}$ such that
\be\label{CONVER}
(\tilde{U}_{1,2;m},\tilde{g}_{1;m})\xrightarrow{C^{1,\beta}} (\mathbb{T}^{2}\times I_{1,2},\tilde{g}_{1}=dx^{2}+\tilde{h}_{1})
\ee   
where, for $x\in I_{1,2}$, $\tilde{h}_{1}(x):=\tilde{h}_{1}|_{\mathbb{T}^{2}\times \{x\}}$ is a $\torus$-symmetric Riemannian metric. 
Note that because the convergence (\ref{CONVER}) is in $C^{1,\beta}$, the ``path" $x\rightarrow \tilde{h}_{1}(x)$ is $C^{1}$. Therefore the second fundamental forms $\tilde{\Theta}_{1}(x):=\tilde{\Theta}_{1}\big|_{\mathbb{T}^{2}\times \{x\}}=\big(\frac{1}{2}\partial_{x} \tilde{h}\big)\big|_{\mathbb{T}^{2}\times \{x\}}$ of the slices $\mathbb{T}^{2}\times \{x\}$ define a continuous ``path" of $\torus$-symmetric, symmetric two-tensors. Denote the mean curvatures by $\tilde{\theta}_{1}(x):=tr_{\tilde{h}_{1}(x)}\tilde{\Theta}_{1}(x)$.
Moreover, also from Proposition \ref{SADP} and Lemma \ref{LL}, there are $C^{1}$-fibrations $f_{m}:U_{1,2;m}\rightarrow I_{1,2}$ such that
\be\label{CONVER2}
\pi_{m}^{-1}(f_{m}^{-1}(x))\xrightarrow{C^{1}} \mathbb{T}^{2}\times \{x\}.
\ee
The $C^{1}$ convergence here is not optimal for the argumentation below as we want to have control on the second fundamental forms of the fibers. However in the technical Proposition \ref{AP2}, which we prove in the Appendix, it is shown that in this situation $f_{m}$ can indeed be chosen to achieve convergence in $C^{2}$ in (\ref{CONVER2}). 
{\it We will assume that this is the case from now on}. 

We want to prove that $\tilde{h}_{1}(x)=\tilde{h}_{0}$. This will follow directly from the next two claims and the identity $\partial_{x}\tilde{h}_{1}(x)=2\tilde{\Theta}_{1}(x)$.

\vs
{\it Claim 1}: If $\tilde{\theta}_{1}(x)=0$ at every slice of $\mathbb{T}^{2}\times I_{1,2}$ then $\tilde{\Theta}_{1}(x)=0$ at every slice of $\mathbb{T}^{2}\times I_{1,2}$. 

\vs
{\it Claim 2}: $\tilde{\theta}_{1}(x)=0$ at every slice of $\mathbb{T}^{2}\times I_{1,2}$.

\vs
\n We prove first {\it Claim 1}. Let $\varphi_{m}:\mathbb{T}^{2}\times I_{1,2}\rightarrow U_{1,2;m}$ be a sequence of diffeomorphisms such that $\varphi_{m}^{*}(\tilde{g}_{1;m})$ converges in $C^{1,\beta}$ to $\tilde{g}_{1}$. Then we can write\footnote{If necessary, $\varphi_{m}$ can be slightly modified to avoid cross terms, as in the metric expression (\ref{PF}).}
\be\label{PF}
\varphi_{m}^{*}(\tilde{g}_{1;m})=\alpha_{m}^{2}dx^{2}+\tilde{h}_{1;m}(x)
\ee
where the real function $\alpha_{m}:\mathbb{T}^{2}\times I_{1,2}\rightarrow \mathbb{R}^{+}$ converges (in $C^{1}$), and as $m\rightarrow \infty$, to the constant function one on $\mathbb{T}^{2}\times I_{1,2}$ and $\tilde{h}_{1;m}$ converges (in $C^{1}$) to $\tilde{h}_{1}$. Let $\tilde{\Theta}_{1;m}(x)$ and $\tilde{\theta}_{1;m}(x)$ be the second fundamental forms and mean curvatures of the slices $\mathbb{T}^{2}\times \{x\}$, as slices in $(\mathbb{T}^{2}\times I_{1,2},\varphi_{m}^{*}(\tilde{g}_{1;m}))$. Then
\be\label{FE}
\partial_{x}\tilde{\theta}_{1;m}=-\Delta_{\tilde{h}_{1;m}}\alpha_{m} +\big(|\tilde{\Theta}_{1;m}|^{2}_{\tilde{h}_{1;m}}+Ric_{\tilde{g}_{1;m}}(\mathfrak{n},\mathfrak{n})\big)\alpha_{m}
\ee
where $\Delta_{\tilde{h}_{1;m}}$ is the $\tilde{h}_{1;m}$-Laplacian on the slices $\mathbb{T}^{2}\times \{x\}$ and ${\mathfrak n}$ is the unit normal field to the slices.  Let $\zeta(x)$ be a $C^{1}$ non-negative real function of one variable with support in $I_{1,2}$ and consider the volume measure on $\mathbb{T}^{2}\times I_{1,2}$, given by $dV_{m}=dA_{\tilde{h}_{1;m}}dx$, where $dA_{\tilde{h}_{1;m}}$ is the area element of $\tilde{h}_{1;m}$ on every slice. Multiplying (\ref{FE}) by $\zeta dV_{m}$ and integrating we obtain: for the integral of the left hand side and after integration by parts in the variable $x$
\be\label{L1}
-\int_{\mathbb{T}^{2}\times I_{1,2}} \big((\partial_{x} \zeta)\tilde{\theta}_{1;m}+\zeta 2\alpha_{m}(\tilde{\theta}_{1;m})^{2}\big)dV_{m} 
\ee
where we used $\partial_{x} dA_{\tilde{h}_{1;m}}/dA_{\tilde{h}_{1;m}}=\alpha_{m}\tilde{\theta}_{m}$, and, for the integral of the first term of the right hand side exactly the value zero because $\zeta$ is constant over every slice. As $\alpha_{m}\xrightarrow{m\rightarrow \infty} 1$ and $\tilde{\theta}_{1;m}\xrightarrow{m\rightarrow \infty} \tilde{\theta}_{1}$ (in $C^{1}$ and $C^{0}$ resp. and all over $\mathbb{T}^{2}\times I_{1,2}$) we conclude that if $\tilde{\theta}_{1}=0$ then (\ref{L1}) goes to zero, and that this is so for any $\zeta$. Therefore the integral of the second term in the right hand side of (\ref{FE}), namely 
\ben
\int_{{\mathbb{T}^{2}\times I_{1,2}}} \zeta \big(|\tilde{\Theta}_{1;m}|^{2}_{\tilde{h}_{1;m}}+Ric_{\tilde{g}_{1;m}}(\mathfrak{n},\mathfrak{n})\big)\alpha_{m}dV_{m}
\een
must go to zero independently of $\zeta$. But  $Ric_{\tilde{g}_{1;m}}(\mathfrak{n},\mathfrak{n})\geq 0$ for every $m$ and thus, in the limit, we must have $\int \zeta |\tilde{\Theta}_{1}|_{\tilde{h}_{0}}dA_{\tilde{h}_{0}}dx=0$ for every $\zeta$. Hence $\tilde{\Theta}_{1}=0$ as claimed.

\vs
We prove now {\it Claim 2}. We show first the impossibility of having, for some $\bar{x}$, $\tilde{\theta}_{1}(\bar{x})<0$. After that we prove the impossibility of having $\tilde{\theta}_{1}(\bar{x})>0$. To do so we will appeal to the following standard fact. Fact 1: {\it Let $S\subset M$ be a hypersurface on a manifold $M$ with a unit-normal field $\mathfrak{n}$. Let $p\in M$ and $\gamma$ a geodesic segment starting at $S$ in the direction of $\mathfrak{n}$, ending at $p$ and with $dist(p,S)=length(\gamma)$. If $\theta|_{S}\leq \theta_{0}<0$ and $Ric\geq0$ all over a neighborhood of $\gamma$, then $length(\gamma)\leq 2/|\theta_{0}|$.}  

\vs
$\bullet$ {\it Suppose that, for some $\bar{x}$, $\tilde{\theta}_{1}(\bar{x})<0$}. Then by (\ref{CONVER2}) we conclude that there is $m_{2}\geq m_{1}$ such that for every $m\geq m_{2}$ we have $\theta_{1;m}|_{f_{m}^{-1}(\bar{x})}<\tilde{\theta}_{1}(\bar{x})/2$, where $\theta_{m}(x)$ is the mean curvature of $f_{m}^{-1}(x)$, namely $\pi^{*}_{m}(\theta_{m})=\tilde{\theta}_{m}$. 
But note that: the solid torus $M(f_{m}^{-1}(\bar{x}))$ lies inside $M(\st_{1;m})$ which is a region of non-negative Ricci, that $\partial M(f_{m}^{-1}(\bar{x}))$ is $f_{m}^{-1}(\bar{x})$ and finally by Step {\bf A} that $rad_{g_{1;m}}(M(f^{-1}_{m}(\bar{x})))\rightarrow \infty$. This easily contradicts Fact 1, as then for any $m\geq m_{2}$ there is a point $p_{m}$ and a geodesic segment in $M(f_{m}^{-1}(\bar{x}))$ starting at $f^{-1}_{m}(\bar{x})$ and ending at $p_{m}$ of $g_{1,m}$-length equal to $rad_{g_{1;m}}(M(f_{m}^{-1}(\bar{x})))$ and therefore realizing the $g_{1;m}$-distance from $p_{m}$ to $f_{m}^{-1}(\bar{x})$.

\vs
$\bullet$ {\it Suppose that, for some $\bar{x}$, $\tilde{\theta}_{1}(\bar{x})>0$}. Again by (\ref{CONVER2}) we conclude that there is $m'_{2}\geq m_{1}$ such that for every $m\geq m'_{2}$ we have $\theta_{1;m}|_{f_{m}^{-1}(\bar{x})}>\tilde{\theta}_{1}(\bar{x})/2$. We will prove that there is a sequence of geodesic segments $\eta_{m}$, for $m\geq m_{3}$, lying entirely inside $\rt\setminus (M(\st_{1;m})^{\circ}\cup \overline{B_{g}(o,r_{0})})$, starting at $\st_{1;m}$ and ending at a point $p_{m}$ and with 
\begin{gather*}
d^{\rt}_{g_{1;m}}(p_{m},\st_{1;m})=length_{g_{1;m}}(\eta_{m}),\\
length_{g_{1;m}}(\eta_{m})\xrightarrow{m\rightarrow \infty} \infty. 
\end{gather*}
About this sequence we make two crucial remarks: first, the geodesic $\eta_{m}$ will lie entirely in the open set $\rt\setminus \overline{B_{g}(o,r_{0})}$ where the Ricci curvature is non-negative; second, the mean curvature at the initial point of $\eta_{m}$ in $\st_{1;m}$, and in the direction $\eta_{m}'$ (which is opposite to the one used to define $\tilde{\theta}_{1}(x)$) is less or equal than $-\tilde{\theta}_{1}(\bar{x})/2<0$. That $\tilde{\theta}_{1}(\bar{x})$ cannot be positive, contrary to what was assumed, will follow directly from these two remarks and Fact 1. We move then to prove the existence of such sequence.    

Recall that {\it a ray} is an infinite-length geodesic diffeomorphic to $[0,\infty)=\mathbb{R}^{+}\cup \{0\}$ minimizing the distance between any two of its points. Let ${\mathfrak R}_{r_{0}}$ be the set of rays $\xi$ in $(\rt,g)$ starting at a {\it base point $b(\xi)$} in $\partial B_{g}(o,r_{0})$ and lying entirely inside the closed set $\rt\setminus B_{g}(o,r_{0})$. The family ${\mathfrak R}_{r_{0}}$ is easily seen to be non-empty and the union of the rays in ${\mathfrak R}_{r_{0}}$ to be a closed set in $\rt$. Moreover observe the following simple fact about ${\mathfrak R}_{r_{0}}$ to be used later. {\it Consider a sequence $\gamma_{j}$ of geodesic segments lying entirely in $\rt\setminus B_{g}(o,r_{0})$, having one of its end points in $\partial B_{g}(o,r_{0})$ and minimizing the distance between its two extreme points. If $length_{g}(\gamma_{j})\rightarrow \infty$, then there is a subsequence of $\gamma_{j}$ converging (on compact sets of $\rt$) to a ray in ${\mathfrak R}_{r_{0}}$.}  

Let ${\mathcal P}_{L}$ be the set of points in the rays of ${\mathfrak R}_{r_{0}}$ lying at a $g$-distance $L$ from the base point of the ray to which they belong, more precisely
\ben
{\mathcal P}_{L}=\{p\in \xi\in {\mathfrak R}_{r_{0}}/d^{\rt}_{g}(p,b(\xi))=L\}
\een
\begin{figure}[h]
\centering
\includegraphics[width=8cm,height=8cm]{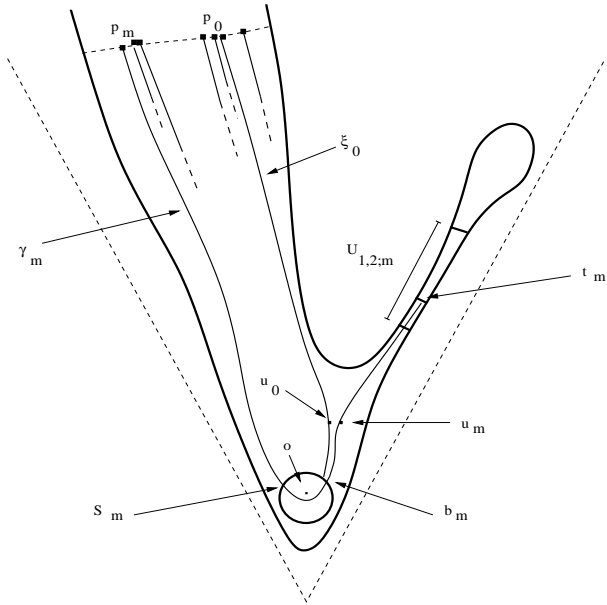}
\caption{Representation of the construction in the proof of Step {\bf B}. In terms of length it more economic to go from $p_{0}$ to $t_{m}$ using the path: $p_{0}\xrightarrow{{\rm by} \xi_{0}} u_{0}\xrightarrow{{\rm by\ short\ curve}} u_{m}\xrightarrow{{\rm by} \gamma_{m}} t_{m}$, rather than going from $p_{0}$ to $t_{m}$ along $\gamma_{m}$. }
\label{CF}
\end{figure} 
Now, for every $m$ there is $L_{m}>0$ sufficiently big with the following properties\footnote{If $M(\st_{1;m})\subset B_{g_{k_{1;m}}}(o,\bar{L}_{m})$ then take $L_{m}=m\bar{L}_{m}10^{k_{1;m}}$.}:
\begin{enumerate}
\item[\bf P1.] ${\mathcal P}_{L_{m}}\subset \big(\rt\setminus M(\st_{1;m})\big)$, and,
\item[\bf P2.] $d^{\rt}_{g_{1;m}}({\mathcal P}_{L_{m}},\st_{1;m})\geq m$ (note that the distance is with respect to $g_{1;m}$).
\end{enumerate} 
Let $\gamma_{m}$ be a geodesic segment from a point $t_{m}$ in $\st_{1;m}$ to a point $p_{m}$ in ${\mathcal P}_{L_{m}}$ and realizing the $g_{1;m}$-distance between the closed sets $\st_{1;m}$ and ${\mathcal P}_{L_{m}}$ of $\rt$. Because of {\bf P1} such segment must lie entirely in $\rt\setminus M(\st_{1;m})^{\circ}$ but it is a priori not evident that it will not intersect $B_{g}(o,r_{0})$. We show now that there is $m_{3}\geq m'_{2}$ such that for every $m\geq m_{3}$ $\gamma_{m}\cap \overline{B_{g}(o,r_{0})}=\emptyset$. With this information and ${\bf P2}$ we can conclude that $\eta_{m}:=\gamma_{m}$ is the sequence we claimed for and the {\it claim} 2 will be finished. Suppose on the contrary that there is a subsequence (denoted again by $\gamma_{m}$) such that $\gamma_{m}\cap \overline{B_{g}(o,r_{0})}\neq\emptyset$.
In this case $\gamma_{m}\cap \overline{B_{g}(o,r_{0})}$, as a closed set in $\gamma_{m}$, has a point $b_{m}$ nearest to $t_{m}$ and a point $s_{m}$ nearest to $p_{m}$.
Let $\hat{\gamma}_{m}$ be the piece of $\gamma_{m}$ enclosed between $t_{m}$ and $b_{m}$. 
Obviously $\hat{\gamma}_{m}$ 
lies inside $\rt\setminus B_{g}(o,r_{0})$. Therefore, as commented above, the sequence $\hat{\gamma}_{m}$ has a subsequence (denoted again by $\hat{\gamma}_{m}$) converging to a ray $\xi_{0}$ (on compact sets of $\rt$). Let $u_{0}$ be a point in $\xi_{0}$ at a $g$-distance $4r_{0}$ from the base point $b(\xi_{0})$ at $\partial B_{g}(o,r_{0})$. Let $u_{m}$ be a sequence of points in $\hat{\gamma}_{m}$ converging to $u_{0}$. Then for every $\epsilon>0$ there is $m(\epsilon)$ such that for any $m\geq m(\epsilon)$ we have 
\begin{gather*}
d^{\rt}_{g}(u_{0},u_{m})\leq \epsilon,\ \text{and } 4r_{0}-\epsilon\leq d^{\rt}_{g}(u_{m}, b_{m})\leq 4r_{0}+\epsilon
\end{gather*}
Let $p_{0}$ be  in $\xi_{0}$ at a $g$-distance $L_{m}$ from $b(\xi_{0})$, which, by definition, is a point in ${\mathcal P}_{L_{m}}$. Then, if $m\geq m(\epsilon)$ we can write
\begin{align*}
d^{\rt}_{g}({\mathcal P}_{L_{m}},\st_{1;m})&\leq d^{\rt}_{g}(p_{0},t_{m})\leq d^{\rt}_{g}(p_{0},u_{0})+d^{\rt}_{g}(u_{0},u_{m})+d^{\rt}_{g}(u_{m},t_{m})\\
&\leq L_{m}-4r_{0} +\epsilon+ d^{\rt}_{g}(u_{m},t_{m})
\end{align*}
On the other hand
\begin{align*}
d^{\rt}_{g}({\mathcal P}_{m},\st_{1;m})&=d^{\rt}_{g}(p_{m},t_{m})\geq d^{\rt}_{g}(p_{m},s_{m})+d^{\rt}_{g}(b_{m},t_{m})\\
&\geq L_{m}-2r_{0}+d^{\rt}_{g}(b_{m},t_{m})\geq L_{m}-2r_{0}+d^{\rt}_{g}(u_{m},t_{m})
\end{align*}
where we used that $d^{\rt}_{g}(p_{m},s_{m})\geq L_{m}-2r_{0}$ which is easily deduced from the fact that, because $p_{m}\in {\mathcal P}_{L_{m}}$, we have $d^{\rt}_{g}(p_{m},\partial B_{g}(o,r_{0}))\leq L_{m}$.  
The two equations before lead readily to the inequality $2r_{0}\leq \epsilon$ which is impossible if one choses for instance $\epsilon=r_{0}$. A representation of the construction can be seen in Figure \ref{CF}. This finishes the proof of {\it Claim 2} and therefore of Step {\bf B}. \et

\vs
{\bf\small Proof of Step C}. We work here with the subsequence of Step {\bf B}, but to simplify notation still use the subindex $m$. On $U_{1,2;m}$ define the $C^{1}$ vector field $W_{m}=\nabla f_{m}/|\nabla f_{m}|^{2}$ and
on $\tilde{U}_{1,2;m}$ define the lifted function $\tilde{f}_{m}=f_{m}\circ \pi_{m}$ and the lifted vector field $\tilde{W}_{m}=\nabla \tilde{f}_{m}/|\nabla \tilde{f}_{m}|^{2}$. $W_{m}$ and $\tilde{W}_{m}$ define flows $\psi_{m}$ and $\tilde{\psi}_{m}$ on $U_{1,2;m}$ and $\tilde{U}_{1,2;m}$ respectively. Because $df_{m}(W_{m})=1$ and $d\tilde{f}_{m}(\tilde{W}_{m})=1$ the flows $\psi_{m}$ and $\tilde{\psi}_{m}$ take fibers into fibers, that is if $x_{1},\ x_{2}\in I$ then
\ben
\psi_{m}(x_{2}-x_{1},-):f^{-1}_{m}(x_{1})\rightarrow f^{-1}_{m}(x_{2}), \text{ and }
\tilde{\psi}_{m}(x_{2}-x_{1},-):\tilde{f}^{-1}_{m}(x_{1})\rightarrow \tilde{f}^{-1}_{m}(x_{2})
\een
Fix $x_{0}\in I_{1,2}$. Let $\tilde{\chi}_{m}:\torus\rightarrow \tilde{f}^{-1}_{m}(x_{0})$ be chosen (to be concrete) such that as $m\rightarrow \infty$ and as $\tilde{f}^{-1}_{m}(x_{0})\xrightarrow{C^{2}}\torus\times \{x_{0}\}$, $\chi_{m}$ converges in $C^{2}$ to the ``identity"  diffeomorphism: $t\in \torus\rightarrow (t,x_{0})\in \torus\times I_{1,2}$. With the help of $\tilde{\chi}_{m}$ and $\tilde{\psi}_{m}$ one can define $C^{2}$-diffeomorphisms 
\ben
\tilde{\varphi}_{m}:\torus\times I_{1,2}\rightarrow \tilde{U}_{1,2;m},\text{ as  }
\tilde{\varphi}_{m}(t,x)=\tilde{\psi}_{m}(x-x_{0},\tilde{\chi}_{m}(t))
\een 
for which we have $\tilde{\varphi}_{m}(\torus\times \{x\})=\tilde{f}^{-1}_{m}(x)$ and $d\tilde{\varphi}_{m}(\partial_{x})=\tilde{W}_{m}$. 
Moreover as $\tilde{W}_{m}$ is perpendicular to the fibers we have the following form of the pull-back metric
\ben
\tilde{\varphi}^{*}_{m}\tilde{g}_{1;m}=\alpha_{m}^{2}dx^{2}+\tilde{h}_{1;m}(x)
\een
where $\alpha_{m}$ and $\tilde{h}_{1;m}(x)$ (may be different from those in Step {\bf B} but we name them the same) converge in $C^{1}$ to the function identically one and $\tilde{h}_{0}$ respectively. We inspect now the behavior of the length of curves on fibers when we translate them along $\partial_{x}$. Let $\tilde{\curve}_{x_{1}}$ be a curve on $\torus\times \{x_{1}\}$ and $\tilde{\curve}_{x}$ be the transported of $\tilde{\curve}_{x_{1}}$ by $\partial_{x}$ to $\torus\times \{x\}$. 
Then, as $\partial_{x} \tilde{h}_{m}=2\alpha_{1,m}\tilde{\Theta}_{1;m}$ we obtain the following direct estimate
\be\label{LCL}
|\partial_{x} length_{\tilde{h}_{1;m}(x)}(\tilde{\curve}_{x})|\leq \frac{\big(\sup_{\torus\times \{x\}}|\tilde{\Theta}_{1;m}|\big)}{2}length_{\tilde{h}_{1;m}(x)}(\tilde{\curve}_{x})
\ee
But, because of Step {\bf B}, $|\tilde{\Theta}_{1;m}|_{\tilde{h}_{1;m}}\xrightarrow{m\rightarrow \infty} 0$ (uniformly on $\torus\times I_{1,2}$) we deduce that: {\it for every $1>\nu>0$ there is $m(\nu)$ such that if $m\geq m(\nu)$ and $x_{1},x_{2}\in I_{1,2}$, then} 
\be\label{LE}
(1-\nu) length_{\tilde{h}_{1;m}(x_{1})}(\tilde{\mathscr C}_{x_{1}})\leq  length_{\tilde{h}_{1;m}(x_{2})}(\tilde{\mathscr C}_{x_{2}})\leq (1+\nu)  length_{\tilde{h}_{1;m}(x_{1})}(\tilde{\mathscr C}_{x_{1}})
\ee
Now, from (\ref{LE}) and noting that the result of transporting a curve $\curve_{x_{1}}\subset f^{-1}_{m}(x_{1})$ (closed or not) by $W_{m}$ to a curve $\curve_{x_{2}}\subset f^{-1}_{m}(x_{2})$ is the same as the result of lifting $\curve_{x_{1}}$ to an (equal length) curve $\tilde{\curve}_{x_{1}}\subset \torus\times \{x_{1}\}$ by means of $\pi_{m}\circ \tilde{\varphi}_{m}$, transport it by $\partial_{x}$ to a curve $\tilde{\curve}_{x_{2}}$, and then push it down to an (equal length) curve $\curve_{x_{2}}\subset f_{m}^{-1}(x_{2})$, we deduce that if $m\geq m(\nu)$ and $x_{1},x_{2}\in I_{1,2}$ then
\be\label{LE2}
(1-\nu) length_{h_{1;m}(x_{1})}(\curve_{x_{1}})\leq  length_{h_{1;m}(x_{2})}(\curve_{x_{2}})\leq (1+\nu)  length_{h_{1;m}(x_{1})}(\curve_{x_{1}})
\ee  

We are ready to prove that there is $m_{2}$ such that if $m\geq m_{2}$ then $\rm\bf \bar{C}1$ and $\rm\bf \bar{C}2$ holds. We prove first $\rm\bf \bar{C}1$ and then $\rm \bf \bar{C}2$. 

\vs
$\bullet$ First, since the $h_{1;m}(x)$-diameters of the fibers $f^{-1}_{m}(x)$, here denoted by $\Gamma_{1;m}(x)$, are realized by the length of geodesic segments (inside the fiber), then we obtain from (\ref{LE2}) 
\be\label{U}
1-\nu\leq \frac{\Gamma_{1;m}(x_{1})}{\Gamma_{1;m}(x_{2})}\leq 1+\nu
\ee
for any $x_{1},x_{2}\in I_{1,2}$ and $m\geq m(\nu)$. Secondly, in exactly the same way that we proved (\ref{DGHD}) in the example of Section \ref{GHDRE} one can prove the following statement: {\it Given $\Lambda_{1}$ there are $\nu_{0}$ and $\Gamma_{0}$ such that for any Riemannian manifold $(V,g_{V})$ with $|Ric_{g_{V}}|\leq \Lambda_{1}$ and with a $\torus$-fibration $f_{V}:V\rightarrow I_{V}$ ($|I_{V}|\geq 1$) for which}
\ben
1-\nu_{0}\leq \frac{\Gamma_{V}(x_{1})}{\Gamma_{V}(x_{2})}\leq 1+\nu_{0},\ x_{1},x_{2}\in I_{V}, \text{ and }\sup_{x\in I_{V}} \Gamma_{V}(x)\leq \Gamma_{0}
\een
\n {\it where $\Gamma_{V}(x)=diam (f^{-1}_{V}(x))$, we have,}

\be\label{FA}
\frac{1}{6}\inf_{x\in I_{V}} \Gamma_{V}(x)\leq dist_{GH}(V,I_{V})\leq \frac{2}{3} \sup_{x\in I_{V}} \Gamma_{V}(x) 
\ee

Now, take $\Lambda_{1}=100\Lambda_{0}$  where $\Lambda_{0}$ is the coefficient that we assumed in the quadratic curvature decay of $g$, that is in $|Ric_{g}|\leq \Lambda_{0}/r^{2}$. Let $\nu_{0}=\nu_{0}(\Lambda_{1})$ and $\Gamma_{0}=\Gamma_{0}(\Lambda_{1})$. Chose $\nu\leq \min\{1/4,\nu_{0}\}$ and $m_{2}\geq m(\nu)$ (as defined above) and sufficiently big that for any $m\geq m_{2}$ we have $\sup_{x\in I_{1,2}}\Gamma_{1;m}(x)\leq \Gamma_{0}$. If as in $\rm\bf\bar{H}1$, $(U_{1;m},g_{1;m})$ is $\epsilon^{*}$-close in the GH-metric to $(I_{1},|\ |)$, then by the first inequality of (\ref{FA}) (applied\footnote{Note that $|Ric_{g_{k_{1;m}}}|\leq 100\Lambda_{0}$ on $U_{1;m}$.} to $V=U_{1;m}$ and $g_{V}=g_{1;m}$) and by (\ref{U}) we have 
\be\label{FAA}
\sup_{x\in I_{2}}\Gamma_{1;m}(x)\leq \frac{6}{1-\nu}\epsilon^{*}
\ee
Hence by (\ref{FAA}), and the second inequality of (\ref{FA}) (applied\footnote{Note that $|Ric_{g_{k_{2;m}}}|\leq 100\Lambda_{0}$ on $U_{2;m}$.} to $V=U_{2;m}$ and $g_{V}=g_{2;m}$) and recalling that $g_{2;m}=\frac{1}{10^{2}}g_{1;m}$ (implying $\Gamma_{2;m}(x)=\Gamma_{1;m}(x)/10$) we obtain
\ben
dist_{GH}((U_{2;m},g_{2;m}),(I_{2},|\ |))\leq \frac{2}{3}\frac{1}{10}\frac{6}{(1-\nu)}\epsilon^{*}\leq \frac{2}{3}\epsilon^{*}
\een
where the last inequality is because $\nu\leq 1/4$. This shows that ${\rm\bf \bar{C}1}$ holds. 

\vs
$\bullet$ Suppose, as in $\rm\bf \bar{H}2$, that there is a closed $\curve_{1;m}\in \st_{1;m}$ for which it is $length_{g_{1;m}}({\mathscr C}_{1;m})\leq \ell^{*}_{m}$. Let $x_{1}$ be the left point of the interval $I_{1}$ and $x_{2}$ the left point of the interval $I_{2}$. Then the curve $\curve_{1;m}$ belongs to the fiber $f_{m}^{-1}(x_{1})$. 
%
%
Let $\curve_{2;m}$ be the transport of $\curve_{1;m}$ by $W_{m}$ to $f^{-1}_{m}(x_{2})=\st_{2;m}$. By (\ref{LE}) we have 
\ben
length_{g_{2;m}}(\curve_{2;m})=\frac{1}{10}length_{g_{1;m}}(\curve_{2;m})\leq \frac{4}{3}\frac{1}{10} length_{g_{1;m}}(\curve_{1;m})\leq \frac{2}{3}\ell^{*}.
\een
\n This shows that ${\rm\bf \bar{C}2}$ holds.\etp

\vs
\vs
We are ready to prove Theorem \ref{MT}.

\vs
\vs
\n {\bf\large Proof of Theorem \ref{MT}:}\label{PMTP} We will work with the annuli decomposition of Section \ref{SAD}.
The key to the proof of Theorem \ref{MT} is to show that if $i\geq i_{0}$, for some $i_{0}\geq 0$, the manifold $M(T^{2o}_{i+1},T^{2o}_{i})$ is a IIB-manifold. Once this is shown, the proof of Theorem \ref{MT} is as follows. Let $\mathbb{S}^{2}_{\bar{r}}=\partial \mathbb{B}^{3}(o,\bar{r})$ be the ``coordinate sphere" of radius $\bar{r}$ in $\rt$ and let $\bar{r}$ be large enough that
$\mathbb{S}^{2}_{\bar{r}}\subset \rt\setminus M(T^{2o}_{i_{0}})$. Then as\footnote{Note from the properties of annuli decompositions that for any sequence $T^{2}_{j}$ of pairwise different tori in ${\mathcal N}$ we have $d^{\rt}_{g}(o,T^{2}_{j})\rightarrow \infty$. In particular $d^{\rt}_{g}(o,T^{2o}_{i})\rightarrow \infty$ as $i\rightarrow \infty$. This justifies equation (\ref{COV}).} 
\be\label{COV}
\rt\setminus M(T^{2o}_{i_{0}})^{\circ}=\bigcup_{i=i_{0}}^{\infty} M(T^{2o}_{i+1},T^{2o}_{i})
\ee
we have, for some $i_{1}>0$, 
\be\label{COVI}
\mathbb{S}^{2}_{\bar{r}}\subset \bigcup_{i=i_{0}}^{i=i_{1}} M(T^{2o}_{i+1},T^{2o}_{i})
\ee
By Proposition \ref{IIBM}, the right hand side of (\ref{COVI}) is a IIB-manifold if every one of its summands is a IIB-manifold. Therefore $\mathbb{S}^{2}_{\bar{r}}$ bounds a ball in $\cup_{i=i_{0}}^{i=i_{1}} M(T^{2o}_{i+1},T^{2o}_{i})$ and so bounds a ball in $\rt\setminus \{o\}$ because $\rt\setminus \{o\}$ contains $\cup_{i=i_{0}}^{i=i_{1}} M(T^{2o}_{i+1},T^{2o}_{i})$. But $\mathbb{S}^{2}_{\bar{r}}$ does not bound a ball in $\rt\setminus \{o\}$ and we reach a contradiction. 

We move then to prove that there is $i_{0}\geq 0$ such that for any $i\geq i_{0}$, $M(T^{2o}_{i+1},T^{2o}_{i})$ is a IIB-manifold. 

Define $i_{0}$ such that, for every piece $U_{k,l}\subset \rt\setminus M(T^{2o}_{i_{0}})^{\circ}$ we have ($\epsilon^{*},\ell^{*},k^{*}$ below are as in Proposition \ref{CRUC})
\begin{enumerate}
\item $k\geq k^{*}$, and,
\item $(U_{k,l},g_{k})$ is either $\epsilon^{*}$-close in the GH-metric to either an interval or a two-orbifold, and,
\item if  $(U_{k,l},g_{k})$ is $\epsilon^{*}$-close to a two orbifold then the $g_{k}$-length of the fibers ${\mathscr C}$ of the Seifert structure is less or equal than $\ell^{*}$, i.e. $length_{g_{k}}({\mathscr C})\leq \ell^{*}$.
\end{enumerate}
We will use such a $i_{0}$ from now on and show that if $i\geq i_{0}$ then $M(T^{2o}_{i+1},T^{2o}_{i})$ is a IIB-manifold.
Some notation now. If a piece $U_{k,l}$ on $\rt\setminus M(T^{2o}_{i_{0}})^{\circ}$ is $\epsilon^{*}$-close to an interval then we say that the piece is of type $I(\epsilon^{*})$ and if it is not and therefore is $\epsilon^{*}$-close to a two-orbifold then we say that the piece is of type $II(\epsilon^{*})$. 

Let $i\geq i_{0}$, 

\vs
$\bullet$ If $M(T^{2o}_{i+1},T^{2o}_{i})$ does not contain a piece of type $I(\epsilon^{*})$ 
then $M(T^{2o}_{i+1},T^{2o}_{i})$ is a union of Seifert manifolds (with Seifert structures coinciding at any intersection) and therefore a Seifert manifold with two boundary components, $T^{2o}_{i+1}$ and $T^{2o}_{i}$. It follows that in this case $M(T^{2o}_{i+1},T^{2o}_{i})$ is a IIB-manifold.       

\vs
$\bullet$ If $M(T^{2o}_{i+1},T^{2o}_{i})$ contains a piece of type $I(\epsilon^{*})$ then we can distinguish two cases, 

\vs
\ \ \ \ \ \ (i) $M(T^{2o}_{i+1},T^{2o}_{i})$ is itself a piece of type $I(\epsilon^{*})$ (in this case the only $U_{k,l}$-piece), therefore diffeomorphic to $\mathbb{T}^{2}\times I$ and thus a IIB-manifold, or, 

\vs
\ \ \ \ \ \ (ii) $M(T^{2o}_{i+1},T^{2o}_{i})$ is not a piece of type $I(\epsilon^{*})$. 

\vs
We discuss case (ii) now and show that $M(T^{2o}_{i+1},T^{2o}_{i})$ is also in this case, a IIB-manifold. First, denote by  ${\mathcal N}_{i+1,i}$ the set of boundary components, other than $T^{2o}_{i+1}$ and $T^{2o}_{i}$, of the $U_{k,l}$-pieces composing $M(T^{2o}_{i+1},T^{2o}_{i})$. Then any torus $T^{2}$ in ${\mathcal N}_{i+1,i}$ is $\ll$ than $T^{2o}_{i+1}$ but not related, in the order $\ll$, to $T^{2o}_{i}$ (otherwise would be one of the $T^{2o}_{i}$'s).  Recall that for any $T^{2}\in {\mathcal N}_{i+1,i}$ one can associate a maximal chain $\{T^{2},T^{2o}_{i+1}\}\rightarrow \{T^{2},T^{2}_{1},\ldots, T^{2}_{n},T^{2o}_{i+1}\}$ (see notation in Section \ref{AD}).  Define $\widehat{\mathcal N}_{i+1,i}$ as the set of tori $T^{2}$ in ${\mathcal N}_{i+1,i}$ such that 
\begin{enumerate}
\item $T^{2}$ is the boundary component of a $U_{k,l}$-piece of type $I(\epsilon^{*})$, and, 
\item none of the tori $T^{2}_{1},\ldots,T^{2}_{n-1}$ in the chain $\{T^{2},T^{2o}_{i+1}\}\rightarrow \{T^{2},T^{2}_{1},\ldots,T^{2}_{n-1},T^{2o}_{i+1}\}$ is the boundary component of a $U_{k,l}$-piece of type $I(\epsilon^{*})$.
\end{enumerate}   
Then, the set of tori $\{T^{2o}_{i+1},T^{2o}_{i}\}\cup \widehat{\mathcal N}_{i+1,i}$ enclose the region
\ben
M(T^{2o}_{i+1},T^{2o}_{i})\setminus \bigg(\bigcup_{T^{2}\in \widehat{\mathcal N}_{i+1,i}} M(T^{2})\bigg)
\een
which is formed by pieces of type $II(\epsilon^{*})$. Therefore it is a Seifert manifold with at least three boundary components (two of them are $T^{2o}_{i+1}$ and $T^{2o}_{i}$) and hence a IIB-manifold. Now, the tori $T^{2}$ in $\widehat{\mathcal N}_{i+1,i}$ are either of type $\st$ or of type $\nst$, namely either $M(T^{2})$ is a solid torus or not (see beginning of Sec. \ref{PMT}). Consider $\st$ in $\widehat{\mathcal N}_{i+1,i}$. Then, $\st$ is the boundary of a $U_{k,l}$-piece of type $II(\epsilon^{*})$ and, because $i\geq i_{0}$ and the definition of $i_{0}$, the fibers $\{\mathscr C\}$ of the Seifert structure of such piece have $g_{k(\st)}$-length less or equal than $\ell^{*}$. In particular the fibers $\{\mathscr C\}$ on $\st$ (which, as closed curves, are non-contractible in $\st$) have $g_{k(\st)}$-length less or equal than $\ell^{*}$. Summarizing, we would have, $k(\st)\geq k^{*}$ (because $i\geq i_{0}$) and
\begin{enumerate}
\item[\bf H1'.] $(U_{k(\st),l(\st)},g_{k(\st)})$ is $\epsilon^{*}$-close in the GH-metric to an interval, and,
\item[\bf H2'.] There is a curve ${\mathscr C}\subset \st$ (indeed anyone of the $\mathscr C$'s) non-contractible in $\st$ but contractible in $M(\st)$ such that $length_{g_{k(\st)}}({\mathscr C})\leq \ell^{*}$,   
\end{enumerate}
Therefore (and crucially), if the fibers $\{\curve\}$ on $\st$ are contractible inside $M(\st)$ then applying Proposition \ref{CRUC} iteratively, we would obtain a consecutive sequence of pieces of type $I(\epsilon^{*})$ extending to infinity, i.e. a $\mathbb{T}^{2}\times \mathbb{R}^{+}$-end, which is not possible because then $M(\st)$ would not be compact\footnote{Alternatively it would imply the existence of an embedded torus (a section of such end) dividing $\rt$ into two unbounded connected components, which is not possible.}. We conclude that for every $\st$ in $\widehat{\mathcal N}_{i+1,i}$, the fibers $\{\mathscr C\}$ are non-contractible in $M(\st)$. Therefore, recalling the comment at the end of Section \ref{SRTIIB}, the manifold
\ben
M(T^{2o}_{i+1},T^{2o}_{i})\setminus \bigg(\bigcup_{\nst\in \widehat{\mathcal N}_{i+1,i}} M(\nst)\bigg)   
\een
is a IIB-manifold (note the union on the right hand side is on $\nst\in \widehat{\mathcal N}_{i+1,i}$). Finally, for every $\nst\in \widehat{\mathcal N}_{i+1,i}$, $M(\nst)$ is a IIB-manifold as was explained in Section \ref{SRTIIB}. Therefore by Proposition \ref{IIBM}, $M(T^{2o}_{i+1},T^{2o}_{i})$ is a IIB-manifold. This finishes the proof of Theorem \ref{MT}.\label{ETMLF}\ep 
 
\section{ Appendix.} 

\subsection{ Remarks on manifolds and convergence.}\label{RMC}

A three-manifold $M$ is $C^{k+1,\beta}$, $k\geq 1$, $0<\beta<1$ if it is a topological manifold provided with an atlas with transition functions in $C^{k+1,\beta}$. A {\it Riemannian} three-manifold $(M,g)$ is $C^{k,\beta}$ if $M$ is $C^{k+1,\beta}$ and the entries of $g$ in every coordinate system of the $C^{k+1,\beta}$ atlas of $M$ are $C^{k,\beta}$ functions. 

A sequence of $C^{k,\beta}$ Riemannian manifolds $(M_{i},g_{i})$ converges in $C^{k,\beta}$ to a $C^{k,\beta}$ Riemannian manifold $(M,g)$ if there are $C^{k+1,\beta}$-diffeomoprhisms $\varphi_{i}:M\rightarrow M_{i}$ such that
the entries of $\varphi_{i}^{*} g_{i}$ in every coordinate system of the atlas of $M$, converge in $C^{1,\beta}$ to the entries of $g$ in the coordinate system.

There are norms that we will use that do not depend on the coordinates. In particular on a $C^{k,\beta}$ Riemannian manifold $(M,g)$ one can define the $C^{k'+1}_{g}$-norm, $k'\leq k$, of functions as usual as
\ben
\|f\|_{C^{k'+1}_{g}}=\sup_{x\in M} \bigg( \sum_{j=0}^{j=k'+1}|\nabla^{(j)} f|(x)\bigg)
\een
where $\nabla^{(j)}$ is the operator resulting from applying $\nabla$ $j$-times. Note that $\nabla^{(j)} f=\nabla^{(j-1)}df$ and that the $C^{k'+1}_{g}$ norm of $f$ involves only derivatives of $g$ up to order $k'$. In particular the space $C^{2}_{g}$ is well defined on a $C^{1,\beta}$ Riemannian manifold. Moreover one easily has the following property: If $(M_{i},g_{i})$ converges in $C^{1,\beta}$ to $(M,g)$ (via diffeomorphisms $\varphi_{i}$) and $f_{i}$ is a sequence of functions in $M_{i}$, then there is $i_{0}$ such that for any $i\geq i_{0}$ we have $\|\varphi_{i}^{*} f_{i}\|_{C^{2}_{g}}\leq 2\|f_{i}\|_{C^{2}_{g_{i}}}$ (here $\varphi^{*}_{i}f_{i}=f_{i}\circ \varphi_{i}$).  

\subsection{ Some technical propositions.}

The following theorem would be standard if we were working in the smooth category. With low regularity there are some points to check. 
\begin{Proposition}\label{AP1} Let $(M,g)$ be a compact $C^{1,\beta}$-Riemannian manifold with boundary. Suppose that 
$\phi: \torus\times M\rightarrow M$
is a continuous and free action by isometries. Then there exists a $C^{2,\beta}$-diffeomorphism 
$\varphi: M\rightarrow \torus\times I$
such that $\varphi^{*}g=dx^{2}+h(x)$ where $h(x)$ is a $C^{1,\beta}$-path of $\torus$-symmetric, and therefore flat metrics in $\torus$.
\end{Proposition}
\n {\bf Proof:} By \cite{MR1503467} (Thm 6, pg. 411), the set of orbits $\torus(p)=\{\phi(t,p),t\in \torus\}$, $p\in M$, is a foliation of $M$ by $C^{1}$-embedded tori. Let $\torus_{1}\neq \torus_{2}$ be two leaves and let $\gamma_{12}$ be a geodesic segment realizing the distance between them and therefore perpendicular to them. As the action is by isometries the set $\{\phi(t,\gamma_{12}),t\in \torus\}$ is a foliation of the region enclosed by $\torus_{1}$ and $\torus_{2}$ by geodesic segments realizing the distance between $\torus_{1}$ and $\torus_{2}$ and perpendicular to them. As in this argumentation the leaves $\torus_{1}$ and $\torus_{2}$ are arbitrary it follows that any inextensible geodesic perpendicular to one leaf is also perpendicular to any other leaf. Let $\gamma(x)$, $x$ the arc-length, be one of such geodesics. Define $\varphi:\torus\times I\rightarrow M$, $|I|=length(\gamma)$, as
$\varphi(t,x)=\phi(t,\gamma(x))$. By \cite{MR1503467} (in particular (D) in pg. 402) the map $\varphi$ is a $C^{1}$ diffeomorphism.
We have $\varphi^{*} g=dx^{2}+h(x)$ where $h(x)$ is a $C^{0}$-path of $\torus$-symmetric metrics in $\torus$.
Let $(y,z)$ be (local) flat coordinates on $\torus$ which together with $x$ form (local and $C^{1}$) coordinates. The standard Laplacian acting on certain functions $f$ at least can be computed in the coordinates $(x,y,z)$ as $\Delta f=[{\rm det}\, h]^{-1/2} (\partial_{i}(g^{ij}[{\rm det}\, h]^{\frac{1}{2}}\ \partial_{j} f))$ (because ${\rm det}\, h$ is just $C^{0}$). Such is the case\footnote{To justify the $\Delta f$ in these cases multiply by a smooth and arbitrary test function of compact support and integrate by parts.} when $f=x,y$ or $\int^{x}[{\rm det}\, h]^{-1/2}dx$. As ${\rm det}\, h={\rm det}\, h\, (x)$, the coordinates $y$ and $z$ are harmonic (and $C^{1}$) and therefore from standard elliptic regularity also $C^{2,\beta}$ in $M$ (recall for this that $M$ is $C^{2,\beta}$ and $g$ is $C^{1,\beta}$). It remains to see the regularity of $x$. Define a new coordinate by $\bar{x}=\int^{x}[{\rm det}\, h]^{-1/2} dx$. Then $\bar{x}$ is harmonic and because is $C^{1}$, by standard elliptic regularity again, it is $C^{2,\beta}$ in $M$. Therefore $(\bar{x},y,z)$ is an harmonic and $C^{2,\beta}$ coordinate system. Hence in these coordinates the metric coefficients $g_{ij}$ are of class $C^{1,\beta}$. Thus $[{\rm det}\, h]^{1/2}$ is of class $C^{1,\beta}$ and because $x(\bar{x})=\int^{\bar{x}} [{\rm det}\, h]^{1/2} d\bar{x}$ we deduce that $x$ is also $C^{2,\beta}$ in $M$.\ep
\begin{Proposition}\label{AP2} Suppose that a sequence $(U_{m},g_{m})$ with $|Ric_{g_{m}}|\leq \Lambda_{0}$ collapses metrically to $(I,|\ |)$. Then, there is a covering subsequence $(\tilde{U}_{m_{j}},\tilde{g}_{m_{j}})$ (with covering maps $\pi_{m_{j}}$) converging in $C^{1,\beta}$ to a $\torus$-symmetric space $(\torus\times I,dx^{2}+\tilde{h}(x))$, and there is a sequence of functions $f_{m_{j}}:U_{m_{j}}\rightarrow \mathbb{R}$, such that ${f}_{m_{j}}\circ \pi_{m_{j}}:\torus\times I\rightarrow \mathbb{R}$ converges in $C^{2}$ to the coordinate function $x$. In particular, fixed a value of $x$, $\pi_{m_{j}}^{-1}(f_{m_{j}}^{-1}(x))$ converges in $C^{2}$ to the slice $\torus\times \{x\}$.
\end{Proposition} 
\n {\bf Proof:} The first part of the claim, i.e. the existence of the covering subsequence is known to us from Lemma \ref{LL}. Thus assume that 
\ben
(\tilde{U}_{m_{j}},\tilde{g}_{m_{j}})\xrightarrow{C^{1,\beta}} (\torus\times I,\tilde{g}=dx^{2}+\tilde{h}(x))
\een
Following \cite{MR767363} (see also \cite{MR1145256} pg. 336), for every $\epsilon>0$ there are\footnote{We remark that this useful smoothing procedure has been used recurrently in \cite{MR950552} as it greatly simplifies the arguments. Our use does no differ much from the purposes it was used there.} smoothings $g^{\epsilon}_{m_{j}}$ of $g_{m_{j}}$ such that 
\be\label{ESII}
dist_{Lip}(g_{m_{j}},g^{\epsilon}_{m_{j}})\leq \epsilon,\ |Ric_{g^{\epsilon}_{m_{j}}}|\leq 2\Lambda_{0},\ |\nabla^{(k)}_{g^{\epsilon}_{m_{j}}} Ric_{g^{\epsilon}_{m_{j}}}|\leq \Lambda_{k}(\epsilon), k\geq 1
\ee
where $dist_{Lip}$ is the Lipschitz distance (see \cite{MR2307192},\cite{MR1145256})\footnote{Note that what makes these estimates useful is that they are independent from the injectivity radius.}. In these smoothed spaces one has the following two properties for fixed $\epsilon$.
\begin{enumerate}
\item[\bf E1.] There is a subsequence of $(\tilde{U}_{m_{j}},\tilde{g}_{m_{j}}^{\epsilon})$ (indexed with $m_{j}$ again but depending on $\epsilon$) converging in $C^{\infty}$ and via diffeomorphisms $\chi_{j}$ to $(\torus\times I, \tilde{g}^{\epsilon}=dx^{2}+\tilde{h}^{\epsilon}(x))$. 
Hence as discussed in Section \ref{RMC} there is $j_{0}(\epsilon)$ such that for every $j\geq j_{0}(\epsilon)$ and sequence of functions $F_{j}$ on $\tilde{U}_{m_{j}}$ we have $\|\chi_{j}^{*}F_{j}\|_{C^{2}_{\tilde{g}^{\epsilon}}}
\leq 2\|F_{j}\|_{C^{2}_{\tilde{g}^{\epsilon}_{m_{j}}}}$.
\item[\bf E2.] From Lemma 1.6 (pg. 336) in \cite{MR984756}, there are fibrations $f_{m_{j}}^{\epsilon}:U_{m_{j}}\rightarrow I$ such that, for all $k\geq 1$, 
$\|f_{m_{j}}^{\epsilon}\circ \pi_{m_{j}}\|_{C^{k}_{\tilde{g}^{\epsilon}_{m_{j}}}}\leq C'_{k}(\epsilon)$.
Moreover $f_{m_{j}}^{\epsilon}\circ \pi_{m_{j}}$ converges in $C^{1}$ to the function $x$ in $(\torus\times I, dx^{2}+\tilde{h}^{\epsilon}(x))$ and, because of the estimates before, the convergence is also in $C^{\infty}$. In particular 
$
\lim \|\chi_{j}^{*}(\pi_{m_{j}}\circ f_{m_{j}})-x\|_{C^{2}_{\tilde{g}^{\epsilon}}}\xrightarrow{j\rightarrow \infty} 0
$
\end{enumerate}
And if we make $\epsilon\rightarrow 0$ we have, because of the first two terms of (\ref{ESII}), the following property.
\begin{enumerate}
\item[\bf E3.] As $\epsilon\rightarrow 0$, the spaces $(\torus\times I,dx^{2}+\tilde{h}^{\epsilon}(x))$ converge in $C^{1,\beta'}$ ($\beta'<\beta$) and via diffeomorphisms $\varphi_{\epsilon}$ to $(\torus\times I, dx^{2}+\tilde{h}(x))$. Moreover by Prop. \ref{AP1} the $C^{2}$-coordinates $x$ in them converge in $C^{2}$ to the (by Prop. \ref{AP1}) $C^{2}$-coordinate $x$ in the limit space.
\end{enumerate}
From {\bf E1} and {\bf E3} we immediately obtain that, for every $\epsilon(i)=1/i$, $i=1,2,3,\ldots$ one can find $m_{j(i)}$ with $j(i)\geq j_{0}(\epsilon(i))$, in such a way that the subsequence $(\tilde{U}_{m_{j(i)}},\tilde{g}^{\epsilon(i)}_{m_{j(i)}})$ converges in $C^{1,\beta'}$ and via the diffeomorphisms $\chi_{j(i)}\circ \varphi_{\epsilon(i)}$ to $(\torus\times I, dx^{2}+\tilde{h}(x))$. Then we have

\begin{align*}
\!\!\!\|\varphi_{\epsilon(i)}^{*}\chi_{j(i)}^{*}(\pi_{m_{j(i)}}\circ  f^{\epsilon(i)}_{m_{j(i)}})- x\|_{C^{2}_{\tilde{g}}}
&\leq \|\varphi_{\epsilon(i)}^{*}\chi_{j(i)}^{*}(\pi_{m_{j(i)}}\circ f^{\epsilon(i)}_{m_{j(i)}})-\varphi_{\epsilon(i)}^{*}x+
\varphi_{\epsilon(i)}^{*}x - x\|_{C^{2}_{\tilde{g}}}\\
&\leq 2\|\chi_{j(i)}^{*}(\pi_{m_{j(i)}}\circ f^{\epsilon(i)}_{m_{j(i)}})-x\|_{C^{2}_{\tilde{g}^{\epsilon(i)}_{m_{j(i)}}}}+\|\varphi_{\epsilon(i)}^{*}x - x\|_{C^{2}_{\tilde{g}}}.
\end{align*}
where the last term tends to zero as $i\rightarrow \infty$. \ep

\vs
\n {\bf Proof of Lemma \ref{LL}:} The result is a straightforward consequence of the assumption that $(M_{i},g_{i})\in {\mathcal M}({\mathfrak N}_{0})$ for some fixed ${\mathfrak N}_{0}$ and the results in \cite{MR950552}. There are however some technical points which is better to clarify and these have to do with the fact that several metrics are involved at the same time. Fukaya's proofs of course will not repeated here and we refer the reader to his articles for full information.

For every ``$i$" let $\breve{M}^{\circ}_{i}(p_{i})$ be the connected component of $M_{i}^{\circ}$ containing $p_{i}$ and let $d_{i}=d_{g_{i}}^{M_{i}}$. We let $M^{\underline{\epsilon}}_{i}(p_{i}):=\breve{M}^{\circ}_{i}(p_{i})\setminus {\mathcal T}_{d_{i}}(\partial M_{i},\underline{\epsilon})$ and similarly for $M^{\overline{\epsilon}}_{i}(p_{i})$. From Proposition \ref{PCOM} we can take a subsequence (index again by ``$i$") such that $(M_{i}^{\underline{\epsilon}}(p_{i}),d_{i})$ converges to a compact metric space $(X^{\underline{\epsilon}},d^{\underline{\epsilon}})$. The subsequence can be chosen in such a way that $M^{\overline{\epsilon}}_{i}(p_{i})$ converges (as a compact set) to $X^{\overline{\epsilon}}\subset X^{\underline{\epsilon}}$. 
We keep using this sequence in the following. 

Following \cite{MR950552} (\footnote{Recall the discussion after the statement of Lemma \ref{LL}.}), for every $x\in X^{\overline{\epsilon}}$ there is $\delta(x)\leq (\overline{\epsilon}-\underline{\epsilon})/2$ such that $(B_{d^{\underline{\epsilon}}}(x,\delta(x)),d^{\underline{\epsilon}})$ is locally isometric to a model space {\bf I.a}, {\bf I.b}, {\bf II.a} or {\bf II.b}. Consider then in $B_{d^{\underline{\epsilon}}}(x,\delta(x))$ the corresponding Riemannian metric and denote it by $g^{\underline{\epsilon}}$. In addition to this information, there is a sequence of points $q_{i}\in M^{\overline{\epsilon}}_{i}(p_{i})$ with $q_{i}\rightarrow x$, such that $(\overline{B_{g_{i}}(q_{i},\delta(x))}, g_{i})$ converges in the GH-topology to $(\overline{B_{d^{\underline{\epsilon}}}(x,\delta(x))},g^{\underline{\epsilon}})$. 

Now, using the compactness of $X^{\overline{\epsilon}}$ one can pick points $x_{1},\ldots,x_{J}$ in $X^{\overline{\epsilon}}$ such that the balls $B_{d^{\underline{\epsilon}}}(x_{j},\delta(x_{j})/4)$, $j=1,\ldots,J$, cover $X^{\overline{\epsilon}}$. Assume that $p_{i}$ converges to a point $x_{0}$, that points $p_{j,i}$ converge to $x_{j}$ and that the union $\bigcup_{j=1}^{j=J} B_{d^{\underline{\epsilon}}}(x_{j},\delta(x_{j})/4)$ is connected (otherwise take the connected component of the union containing $x_{0}$). Then, it is direct to check that
\ben
(\bigcup_{j=1}^{j=J}\overline{B_{g_{i}}(p_{j,i},\delta(x_{j}))}, g_{i})\xrightarrow{\rm GH} (\bigcup_{j=1}^{j=J} \overline{B_{d^{\underline{\epsilon}}}(x_{j},\delta(x_{j}))}, g^{\underline{\epsilon}}). 
\een
Then, one can use the local construction in \cite{MR950552} (pg. 19, based on \cite{MR873459}) to find $C^{1}$ functions
\ben
f_{i}: \bigcup_{j=1}^{j=J}\overline{B_{g_{i}}(p_{j,i},3\delta(x_{j})/4)}\rightarrow X^{\underline{\epsilon}}
\een
(but non-surjective) satisfying  the properties in Theorem 0.12 of  \cite{MR950552} (pg. 3) and with range covering 
$\bigcup_{j=1}^{j=J} \overline{B_{d^{\underline{\epsilon}}}(x_{j},\delta(x_{j})/2)}$. The $\Omega_{i}$'s and the space $(X,d)$ are finally defined as
\ben
\Omega_{i}:=f_{i}^{-1}(\bigcup_{j=1}^{j=J} \overline{B_{d^{\underline{\epsilon}}}(x_{j},\delta(x_{j})/2)})\ \text{ and as }
(X,d):=(\bigcup_{j=1}^{j=J} \overline{B_{d^{\underline{\epsilon}}}(x_{j},\delta(x_{j})/2)},g^{\underline{\epsilon}})
\een
with $(X,d)$ satisfying ${\bf D1}$ and ${\bf D2}$ by construction. We then have $(\Omega_{i},g_{i})\xrightarrow{\rm GH} (X,d)$ and $f_{i}:(\Omega_{i}, g_{i})\rightarrow (X,d)$ with the properties ${\bf I}$ which correspond in our case to properties (0.13.1) and (0.13.2) of Theorem 0.12 of  \cite{MR950552}. 

We discuss now how to show {\bf II} in {\it case {\bf D1}}. The {\it case} {\bf D2} is done along similar lines an as we will not use it in this article the proof is left to the reader. Take covers $(\tilde{\Omega}_{i},\tilde{g}_{i})$ to have the injectivity radius at one point controlled away from zero. Leave aside for a moment the issue of the existence of such cover. As $|Ric_{\tilde{g}_{i}}|\leq\Lambda_{0}$ we can take a convergent subsequence, say to $(\tilde{\Omega},\tilde{g})$. The group of Deck-covering transformations of $\tilde{\Omega}_{i}$ converge necessarily to a closed group $G$ of isometries of the limit space  $(\tilde{\Omega},\tilde{g})$. On the other hand, for any $x\in X\setminus Sing(X)$, the fiber $\pi_{i}^{-1}(f^{-1}_{i}(x))$ which covers the torus $f^{-1}_{i}(x)$ under $\pi_{i}$, converges to a torus, say $\tilde{T}^{2}(x)\subset \tilde{\Omega}$. The group $G$ acts effectively\footnote{The action is effective because if an isometry of $G$ leaves every point of $\tilde{T}^{2}(x)$ invariant then it must be the identity as an isometry in $\tilde{\Omega}$.}
by isometries on $\tilde{T}^{2}(x)$ and its quotient is a point. It follows that $G$ is a torus. 

To show that there are covers as mentioned before, observe that, from Lemma \ref{LL}, any ``sufficiently collapsed" manifold (of bounded diameter and curvature) must possess at least one small and non-contractible loop. Now, in {\it case} {\bf D1}, the manifolds $\Omega_{i}$ are diffeomorphic to either $\torus\times \mathbb{I}$ or $\mathbb{B}^{2}\times \mathbb{S}^{1}$ whose fundamental groups are $\mathbb{Z}\times \mathbb{Z}$ and $\mathbb{Z}$ respectively.
In either case one can then take (controlled) covers having no non-contractible and small loops. In this way the cover is necessarily non-collapsed. 
%
\ep

}

\bibliographystyle{plain}

\bibliography{Master.bib}

\end{document}